\documentclass[11pt]{article}
	
	\newcommand{\blind}{0}
	
    \makeatletter
    \renewcommand\section{\@startsection {section}{1}{\z@}%
                                       {-3.5ex \@plus -1ex \@minus -.2ex}%
                                       {2.3ex \@plus.2ex}%
                                       {\normalfont\fontfamily{phv}\fontsize{16}{19}\bfseries}}
    \renewcommand\subsection{\@startsection{subsection}{2}{\z@}%
                                         {-3.25ex\@plus -1ex \@minus -.2ex}%
                                         {1.5ex \@plus .2ex}%
                                         {\normalfont\fontfamily{phv}\fontsize{14}{17}\bfseries}}
    \renewcommand\subsubsection{\@startsection{subsubsection}{3}{\z@}%
                                        {-3.25ex\@plus -1ex \@minus -.2ex}%
                                         {1.5ex \@plus .2ex}%
                                         {\normalfont\normalsize\fontfamily{phv}\fontsize{14}{17}\selectfont}}
    \makeatother
	
	\usepackage{amsmath}
	\usepackage{amssymb}
	\usepackage{graphicx}
	\usepackage{enumerate}
	\usepackage{natbib} 
	\usepackage{url} 
	\usepackage{algpseudocode}
    \usepackage{algorithm2e}
    \usepackage{arevmath}     
    
	\usepackage{caption}
    \usepackage{subcaption}
     \usepackage{multirow}

\begin{document}
    \def\spacingset#1{\renewcommand{\baselinestretch}%
			{#1}\small\normalsize} \spacingset{1} 

		\if0\blind
		{
			\title{\LARGE\bf A Two Stage Stochastic Optimization Model for Port Infrastructure Planning}
			\author{Sanjeev Bhurtyal$^a$, Sarah Hernandez$^a$, Sandra Eksioglu$^b$, Manzi Yves$^a$\\  
			$^a$ Department of Civil Engineering, University of Arkansas, Fayetteville, AR \\
			$^b$ Department of Industrial Engineering, University of Arkansas, Fayetteville, AR}
			\date{}
			\maketitle
		} \fi
		
		\if1\blind
		{
       \title{\bf \emph{IISE Transactions} \LaTeX \ Template}
			\author{Author information is purposely removed for double-blind review}
			
\bigskip
			\bigskip
			\bigskip
			\begin{center}
				{\LARGE\bf A Two Stage Stochastic Optimization Model for Port Infrastructure Planning}
			\end{center}
			\medskip
		} \fi
		\bigskip
		
%
\date{}
\maketitle

\begin{abstract}
This paper investigates inland port infrastructure investment planning under uncertain commodity demand conditions. 
A two-stage stochastic optimization is developed to model the impact of demand uncertainty on infrastructure planning 
and transportation decisions. The two-stage stochastic model minimizes the total expected costs, including the capacity 
expansion investment costs associated with handling equipment and storage, and the expected transportation costs. 
To solve the problem, an accelerated Benders decomposition algorithm is implemented. The Arkansas section of the 
McCllean-Kerr Arkansas River Navigation System (MKARNS) is used as a testing ground for the model. Results show that 
commodity
volume and, as expected, the percent of that volume that moves via waterways (in ton-miles) increases with increasing 
investment in port infrastructure. The model is able to identify a cluster of ports that should receive investment 
in port capacity under any investment scenario. The use of a stochastic approach is justified by calculating the 
value of the stochastic solution (VSS).
\end{abstract}

\section{Introduction}
Barge shipments via inland waterways are one of most efficient modes of freight transportation. A standard 15-barge tow can move about 22,500 tons of commodities, which is equivalent to 225 rail cars or 870 tractor-trailer trucks. Similarly, it takes a gallon of fuel to ship one ton of cargo 59 miles by truck, 202 miles by rail, and 514 miles by barge \cite{WaterwaysValue2000}. Hence, freight shipped by waterways tend to be more fuel efficient and have less impacts on the environmental compared to other mode of transportation \cite{WaterwaysValue2000}. Leveraging inland waterways for freight shipments could mitigate the negative impacts associated with increasing commodity flow in the US by alleviating landside bottlenecks.
\par
Inland waterway ports are a critical part of the national freight network, since they carry about half of the domestic waterborne freight in the United States \cite{USACEWaterborne2017}. They contribute to the economy by moving commodities that are key to the Gross Domestic Product (GDP) of many states, such as agriculture for Arkansas \cite{ArkansasAgricultureProfile}. There are both public and private port terminals in US inland waterway and commodity flow data collected by those inland waterway port operators are proprietary \cite{PrivateandPublicTerminals} which makes it is challenging to understand ports' commodity throughput and capacity. Such data are valuable for infrastructure investment planning, such as adding capacity to port servicing roads, waterway dredging schedules, lock and dam maintenance programming. Furthermore, port level data enables port agencies and authorities to identify ports that needs capacity expansion in terms of operational and storage infrastructures. Publicly available statistics are published by the United States (US) Bureau of Transportation Statistics (BTS) via the Port Performance Freight Statistics. The program provides performance statistics for top 25 ports based on overall cargo tonnage, 20-foot equivalent unit (TEU) of container cargo, and dry bulk cargo tonnage  \cite{HuPortPerformance2021}. However, statistics are limited to the top 25 ports, and only few inland waterway ports are included in the list. Additional data for gathering port commodity flows are the US Census Bureau's US Port Data \cite{USPortData}, Commodity Flow Survey (CFS) \cite{Commodityflowsurvey} and Waterborne Commerce Statistics from the US Army Corps of Engineers(USACE) \cite{WaterborneCommerce}. More recently, the USACE has made available monthly commodity flow data at lock locations through the Lock Performance Monitoring System (LPMS) \cite{LPMS}. However, data from above mentioned sources are constrained in their spatial and/or temporal disaggregation which makes it challenging to collect port specific data for inland waterways. For example, the US CFS provides a periodic five-year snapshot of commodity flows through a survey that is expanded to represent the population using statistical sampling procedures. Notably, the LPMS data give monthly tonnage reports for locks but not ports, and many ports may be located between a pair of locks, thus limiting the insights derived for single ports.
\par

While collecting data at locks and dams may address inland waterway performance related questions, it alone does not help answering port specific investment questions. For example, data on commodity tonnage through a single lock may illuminate travel time delays across a given section of waterway. However, it does not provide a means to guide strategic decisions regarding how a limited monetary investment (perhaps less than the amount needed to improve a single lock’s performance) can be best used to alleviate delays along the inland waterway system. In addition, data on commodity tonnage through a section of waterway between a pair of locks provide historical patterns of commodity flows for a series of ports (between the lock pair) but does not provide insight into the ports’ throughput to capacity ratio.  Therefore, it would be difficult to strategically guide capacity expansion investments without additional data on existing port-specific capacity or operational characteristics. An intuitive investment decision using the currently available data may be to expand the capacity of the largest port among a series of ports between two locks. However, such a decision may lead to a local optimal solution such as adding a new shipping berth to a port to alleviate loading. Form a systems point of view, investments on other ports that have greater access to highways and railways, for example, may lead to global optimal solutions such as a wide investment across many ports. Expanding port infrastructure requires large capital investment and tends to target long lifespans, e.g., 25 years or more. Therefore, it is imperative that port capacity expansion investments are neither under nor over invested. Another facet of strategic decision-making in this context regards the potential for growth or decline in port usage as a function of overall freight shipment demand. Thus, decisions about investments in port infrastructure, as any transportation investment, should be evaluated for different scenarios of freight demand that reflect the unknown nature of economically driven trends seen for freight transport. The demand for commodities for waterway transportation is not known in advance, and since demand affects infrastructure expansion decisions, different demand scenarios must be taken into account when making investment decisions. This calls for an investment model that considers several scenarios and provides optimal inland waterway port infrastructure investment solutions when uncertainty is present.
 
The objective of this research is to develop a tool to guide strategic decisions about investments in waterway transportation infrastructure. These decisions are subject to uncertainty observed in the demand for different commodities shipped via waterways, railways and highways; as well as,   transportation costs and transportation network structure. To accomplish the goal, this study proposes a two-stage stochastic optimization (2-SOP) model that minimizes the total of port infrastructure investment cost and the expected commodity transportation cost. The 2-SOP model has two sets of decision variables, the first and second stage variables. First stage decisions determine equipment and storage investments to increase capacity. These decisions are made prior to the realization of the stochastic parameter, demand for commodities. After the realization of uncertainty, second-stage decisions determine the flow of commodities via different modes of transportation \cite{LiuStochasticDef2018}. The use of 2-SOP model, rather than using a deterministic model, avoids underestimation or overestimation of the supply chain costs when stochastic parameters exist. In addition, 2-SOP problem  mimics reality in which decisions such as infrastructure investments are made first and without full information about
uncertainty.
    
\par
We applied the proposed model to the Arkansas section of the Mcclellan-Kerr Arkansas River Navigation System (MKARNS). Several freight demand scenarios are developed to estimate the stochastic demand parameter. The corresponding deterministic equivalent formulation of the stochastic model is a large-scale mixed-integer linear program that is solved using the L-shaped method. The L-shaped method is also known as the stochastic Benders Decomposition. The addition of Knapsack inequalities and Pareto-optimal cuts improved the running time of the algorithm and the quality of solutions found. The remainder of this paper is structured as follows. A review of prior research related to inland waterway port infrastructure is presented followed by the formulation of a two-stage stochastic optimization model. The case study section presents the application of the proposed model to the MKARNS. This is followed by a discussion of the results. The final section provides the main contributions, conclusions, and limitations of the paper along with the potential future work.

\section{Literature Review}
This section presents the literature pertaining to inland waterway infrastructure planning and two-stage stochastic optimization models for transportation planning. In addition, this section also discusses investment models developed for inland waterway infrastructure, highlighting research gaps.

\subsection{Inland Waterway Infrastructure}
\par

Research on the investment planning of inland waterways focuses on the allocation and scheduling of inland waterway infrastructure maintenance. \citep{MahmoudzadehBudget2021} proposed a two-stage deterministic model that allocates funds for channel dredging and lock maintenance by minimizing commodity transportation cost required to serve the demand. In a similar context, \citep{mitchell2013selection} formulated integer programming models to strategically fund dredging projects considering the interdependencies of the project. The solution of the model was shown to increase commodity throughput. \citep{ratick1996risk} developed a risk-based spatial decision support system to address trade-offs in cost and reliability. The model schedules the dredging in such a way that channel reliability is maximized in any time period. These studies focus on channel and lock maintenance as a means to increase throughput; however, the increased throughput must be processed at ports. If infrastructure within a port (i.e. cranes, lifts, elevators, berths, etc.) is not adequate to handle the increased throughput, then handling delays may be introduced, which in turn increases delays at locks. Therefore, to compliment and improve the existing suite of inland waterway planning models, there is a need for decision models to determine optimal port infrastructure investment.

\par
Along the line of port-specific infrastructure investment, the work by \citep{LagoudisAsianPort2014} proposes an investment decision-making process for port infrastructure investment considering uncertainty in commodity throughput. The model identifies the most profitable infrastructure in operation in a port and ranks investment strategies. This study focused on the profitability rather than the serviceability
of a port. Higher profit does not necessarily mean better performance of the transportation system. Therefore, there is a need to focus on optimizing the performance of the whole system.
 \citep{KohPriority2001} developed a heuristic algorithm that determines optimal locations, size, and timing of inland container port developments. The infrastructure needed by ports that handle containers differs (cranes) from that of ports that handle dry bulk (grain elevators). Therefore, the models developed for inland container ports may not be applicable to ports along inland waterways. 
 \citep{whitman2019multicriteria} integrated a dynamic risk-based interdependency
model with weighted multi-criteria decision analysis techniques to rank investment allocation strategies based on economic loss due to disruptive events. The model determines which commodity to allocate resources for but does not get to the level of detail that recommends the type of infrastructure investment that would be needed to handle that commodity. Further, since commodity movements can share infrastructure, investment strategies based solely on commodity movements will likely differ when considering how the infrastructure is shared across commodities.
\par

\subsection{Two-stage stochastic optimization models for transportation planning}
Two-stage stochastic optimization models are used in a wide array of transportation planning problems, including disaster response \cite{barbarosoglu2004two}, supply chain management \cite{MarufuzzamanEksiogluHuang2014TwoStage}, pavement maintenance and rehabilitation  \cite{ameri2019two}, and $CO_{2}$ disposal planning \cite{han2012developing}. Capacity, supply, demand, technology advancement, and budget represent some of the widely used stochastic parameters in such studies. Along the same line of research, two-stage stochastic models are used to solve waterway transportation-related problems such as those related to empty container management \cite{cheung1998two,mittal2013determining}, storage capacity optimization \cite{liu2020two}, yard crane scheduling \cite{zheng2019two}, port management \cite{Aghalari2020b,Aghalari2020a}, and biomass supply chain \cite{marufuzzaman2017designing}. We will focus our discussion on studies related to inland waterway ports.\par

\citep{Aghalari2020a} provided a two-stage mixed-integer linear programming model to optimize towboats and barges assignment in the inland waterway system. This model assumes stochastic supply of commodities and fluctuation of water levels. The model minimizes the sum of the expected cost of using towboats and barges, and the cost of inventory, transportation, processing, and shortage of commodity. Similarly, \citep{Aghalari2020b} proposed a two-stage stochastic programming model to solve the problem of inland waterway barge and towboat management considering perishable commodities with the aim of minimizing the cost of using towboats and barges and the transportation cost. \citep{mittal2013determining} proposed a two-stage stochastic model to determine an optimal number of inland depots that should be opened within a 10-year time frame while maintaining a minimum system cost of repositioning empty containers while satisfying supply and demand requirements. \citep{marufuzzaman2017designing} developed a dynamic multimodal transportation network design that copes with fluctuating biomass supply. The studies by \citet{Aghalari2020a, Aghalari2020b} focused on port management, while the studies by \citep{marufuzzaman2017designing} and \citep{mittal2013determining} focused on the identification of biorefineries and inland container depot location, respectively. To our knowledge, there are no two-stage stochastic models for inland port-specific infrastructure investment planning. Two stage stochastic model is congruous with inland waterway infrastructure planning where infrastructure investment decisions are made prior to the realization of uncertainty. Hence, such models are deemed appropriate for port infrastructure planning. Likewise, since uncertainty in volume of commodity is a deciding factor on port capacity expansion decision, it is vital to consider it in the investment decision models.\par

\subsection{Our contribution to the existing literature}
Our work is distinct from the studies mentioned above, albeit being connected. First, this study fills a critical gap in the literature by considering multi-modal, multi-period, and multi-commodity nature of inland waterway simultaneously to develop a stochastic port-specific investment decision model. Second, this study captures the uncertainty in the volume of commodities moved via the waterways, which is critical for port infrastructure investment planning. While few studies have attempted to develop stochastic model for port management and scheduling, no prior attempt has been found for multi-modal, multi-period, and multi-commodity stochastic model for port-specific infrastructure model.

\section{Problem Description and Mathematical Model Formulation}
The objective of our study is to optimize investment decisions under uncertainty for inland waterway port capacity expansion infrastructure projects. In this paper, we consider the port capacity to include physical space such as the land area for inventory storage and operational equipment such as cranes and grain elevators. The model determines port infrastructure investments and commodity flow in such a way that the total of the investment cost and the expected transportation cost are minimized.
 \par

\begin{figure}[!htb]
    \centering
    \includegraphics[width = 1\linewidth]{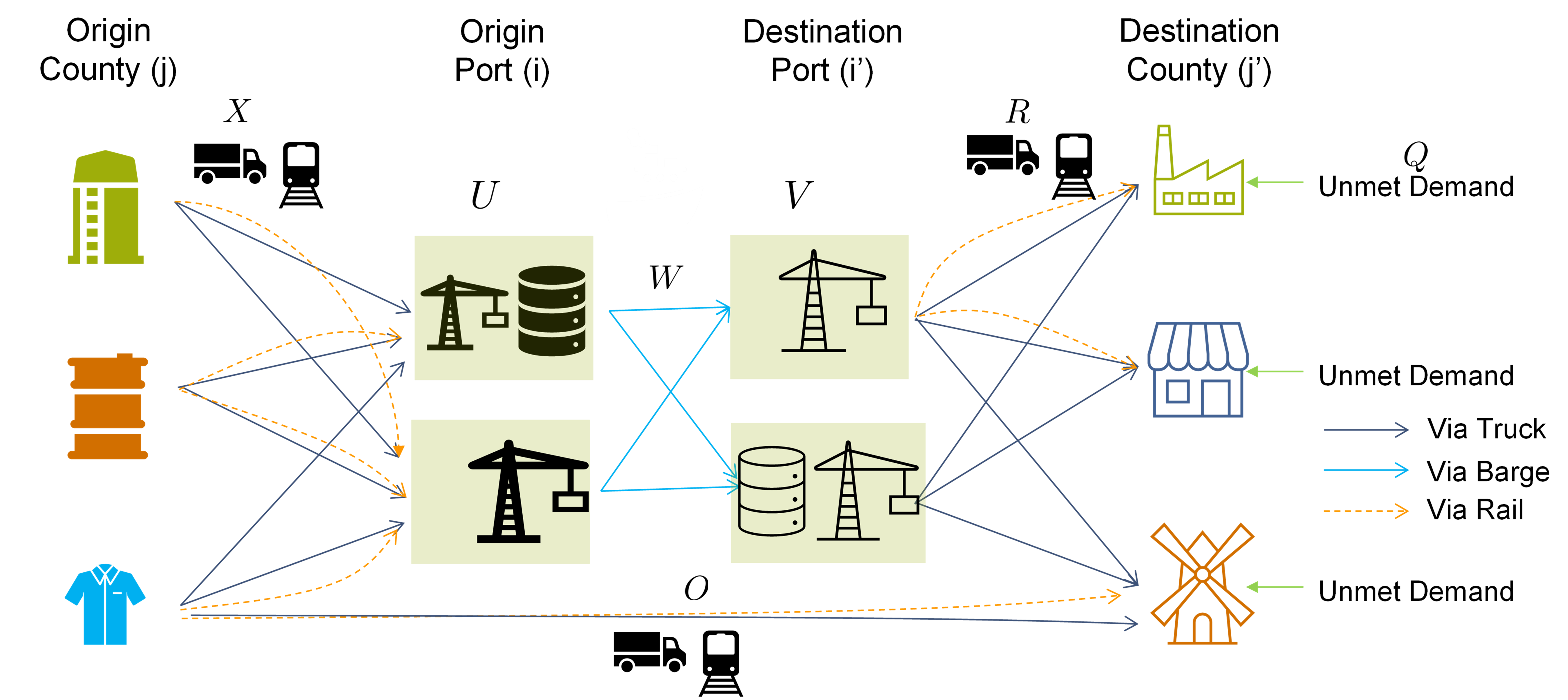}
    \caption{Network representation of supply chain}
    \label{fig:network}
\end{figure}

We consider a supply chain network, where the nodes represent counties that have a positive demand and/or supply for different commodities and ports; and the arcs represent the transportation paths for delivering commodities. Let $W = (N,A)$ denote this transportation network where $N$ is the set of nodes and $A$ is the set of arcs that connect the nodes of the network. The set $N = J \cup I$ ($J \cap I = \emptyset$) represents the counties $J = \{1, 2, \ldots, |J|\}$, from where commodities are shipped and received via a set of land-side transportation modes $T = \{$Truck, Rail$\}$ and a set of ports $I = \{1, 2, \ldots, |I|\}$, from where commodities are transported via barges. The set $C = \{1, 2, 3, \ldots, |C|\}$ represents the set of commodities transported along the network $W$ over a set of time period $P = \{1, 2, \ldots, |P|\}$. We consider two transportation modes between each origin and destination pair. The first path uses rail, truck, and waterways to deliver commodities. Rail and truck are used to deliver commodities to and from each port. Waterways are used to deliver commodities between ports. The second path uses rail and trucks to deliver commodities to and from counties. As an example, consider a supply chain network consisting of three counties that supply commodities via rail, truck and waterways (Figure \ref{fig:network}). Two ports receive the commodities from the counties and ship them via barges to two destination ports. Finally at the destination port, the commodity is shipped to the destination county where there is corresponding demand for that commodity. Any unmet demand is satisfied from other states. \par
\begin{figure}[!htb]
     \centering
     \includegraphics[width = 0.6\linewidth]{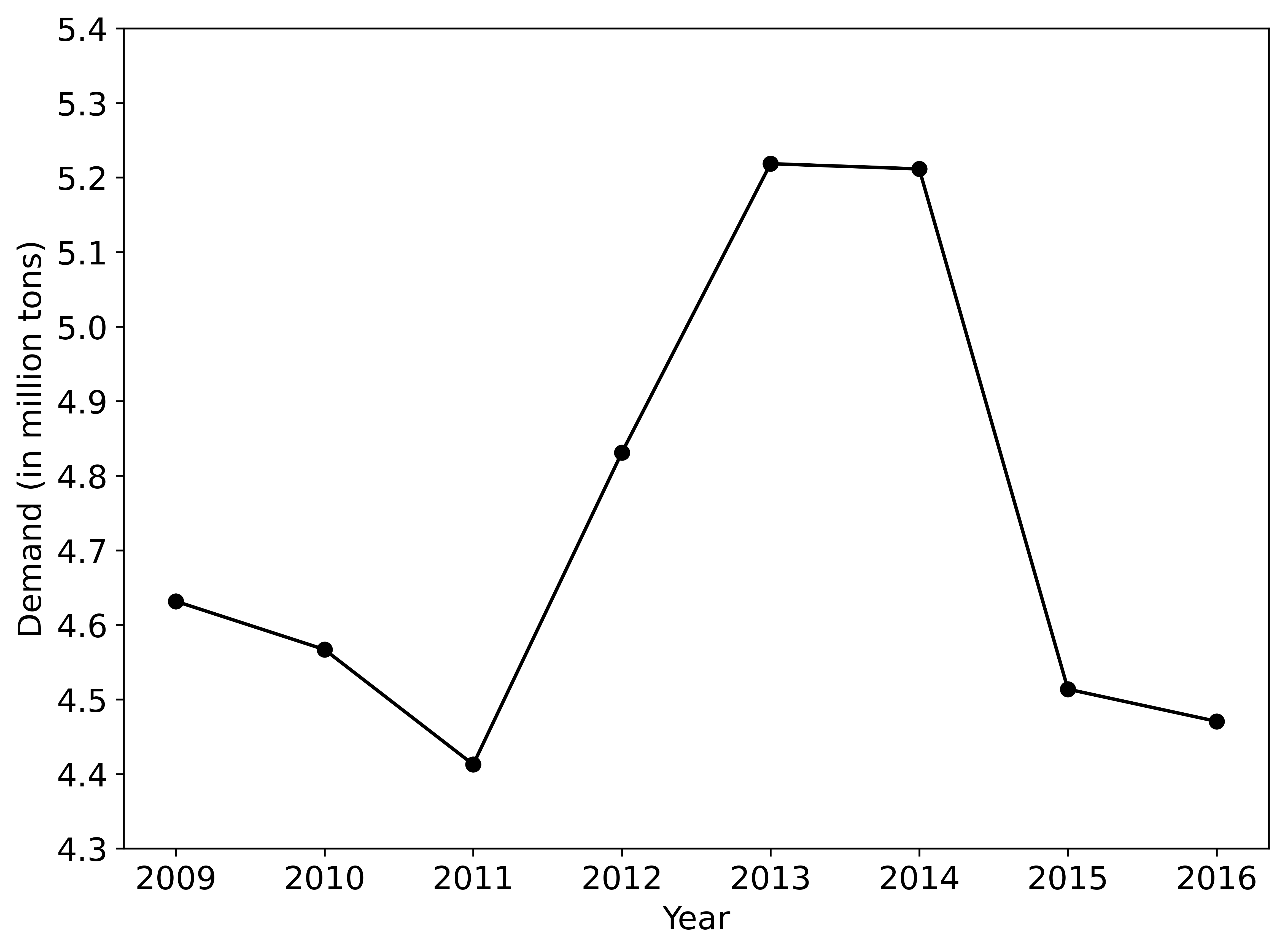}
     \caption{Commodity demand via waterways in Arkansas section of MKARNS}
     \label{fig:stochasticdemand}
 \end{figure}

  
Figure \ref{fig:stochasticdemand} represents commodity demand along the Arkansas section of MKARKNS during 2009 to 2016. This demand changes from year to year, and there is no clear trend. Based on this non-linear and fluctuating demand variation over time, the demand for each commodity is considered stochastic in our model. 

The notations used in the model formulation are summarized as follows:\\
{\underline{Sets:}}
\begin{itemize}
\item $I^{\prime}, I^{{\prime}{\prime}}$  ports of   origins/destinations of shipments ($I = I^{\prime} \cup I^{{\prime}{\prime}})$	
\item $J^{\prime}, J^{{\prime}{\prime}}$ counties of origin / destinations of shipments ($J = J^{\prime} \cup J^{{\prime}{\prime}})$	
	\item $P$ periods in the planning horizon, $p\in P$
	\item $C$ commodities, $c\in C$
	\item $E$ equipment, $e \in E$
	\item $F$ storage facilities, $f \in F$
	\item $S$ scenarios, $s\in S$
\end{itemize}
{\underline{Problem Parameters:}}
\begin{itemize}
    \item $\omega_{s}$ is probability of occurrence of scenario $s$
   \item $\kappa_{e}$ is the cost of equipment $e$
    \item $\iota_{f}$ is the cost of storage facility $f$
    \item $\alpha^t_{ij}$ is the unit cost of transporting commodity  via truck between county $j \in J^{\prime}$ and port $i \in I^{\prime}$ (in \$/ton)
    \item $\alpha^r_{ij}$ is the unit cost of transporting commodity  via rail between county $j \in J^{\prime}$ and port $i \in I^{\prime}$ (in \$/ton)  
    \item $a_{ik}$ is the unit cost of transporting commodity via barge from port $i \in I^{\prime}$ to ${k \in I^{{\prime}{\prime}}}$ (in \$/ton) 
    \item $h_{c}$ is the unit inventory holding cost of commodity $c$
    \item $\mu$ is the unit penalty cost for unmet demand (in \$/ton)
    \item $l^t_{jm}$ is the unit cost for transporting commodity via truck between county $j \in J^{\prime}$ and $m \in J^{{\prime}{\prime}}$ (in \$/ton)
    \item $l^r_{jm}$ is the unit cost for transporting commodity via rail between county $j\in J^{\prime}$ and $m \in J^{{\prime}{\prime}}$ (in \$/ton)
    \item $\beta_{ec}$ takes the value 1 if equipment $e$ can process commodity $c$, and takes the value 0 otherwise
    \item $\Gamma_{fc}$ takes the value 1 if commodity $c$ can be stored in storage facility $f$, and takes the value 0 otherwise
    \item $\delta_i$ takes the value 1 if port $i \in I$ has railway access, and takes the value 0 otherwise
    \item $\gamma_j$ takes the value 1 if county $j \in J$ has railway access, and takes the value 0 otherwise
    \item $B_{fe}$ is the minimum ratio of number of storage facility $f$ to number of equipment $e$
    \item $q_{jcps}$ is the supply availability of commodity $c$ in county $j \in J^{{\prime}}$ in month $p$ in scenario $s$ (in tons)
    \item $d_{jcps}$ is the demand of commodity $c$ in county $j \in J^{{\prime}{\prime}}$ in month $p$ in scenario $s$ (in tons)
    \item $\zeta_{fc}$ is the normalized tonnage of commodity $c$ for inventory in storage facility $f$
    \item $\Lambda_{ec}$ is the normalized tonnage of commodity $c$ for processing in equipment $e$
    \item $l_{f}$ is the storage capacity of storage facility $f$
    \item $k_{if}$ is the existing number of storage facility $f$ at port $i \in I$
    \item $m_{e}$ is the processing capacity of equipment $e$
    \item $n_{ie}$ is the existing number of equipment $e$ at port $i \in I$
\end{itemize}
{\underline{Decision Variables:}}
\begin{itemize}
	\item $Y_{if}$ is the number of storage facility $f$ installed at port $i \in I$
	\item $Z_{ie}$ is the number of equipment $e$ installed at port $i \in I$
   \item $X^t_{jicpes}$ is the tonnage of commodity $c$ shipped via truck from county $j \in J^{\prime}$ to port $i \in I^{\prime}$ and processed using equipment $e$ during period $p$ in scenario $s$
    \item $X^r_{jicpes}$ is the tonnage of commodity $c$ shipped via rail from county $j \in J^{\prime}$ to port $i \in I^{\prime}$ and processed using equipment $e$ during period $p$ in scenario $s$
    \item $U_{icpfs}$ is the tonnage of inventory of commodity $c$ in origin port $i \in I^{\prime}$ in storage facility $f$ in month $p$ in scenario $s$
    \item $W_{ikcpes}$ is the tonnage of commodity $c$ processed using equipment $e$ and transported by barge from ports $i \in I^{\prime}$ and $k \in I^{{\prime}{\prime}}$ in period $p$
    \item $V_{icpfs}$ is the tonnage of inventory of commodity $c$ in destination port $i \in I^{{\prime}{\prime}}$ in storage facility $f$ in month $p$ in scenario $s$
    \item $R^t_{ijcpes}$ is the tonnage of commodity $c$ processed using equipment $e$ and transported from port $i \in I^{{\prime}{\prime}}$ to county $j \in J^{{\prime}{\prime}}$ via truck in month $p$ in scenario $s$
    \item $R^r_{ijcpes}$ is the tonnage of commodity $c$ processed using equipment $e$ and transported from port $i \in I^{{\prime}{\prime}}$ to county $j \in J^{{\prime}{\prime}}$ via rail in month $p$ in scenario $s$
    \item $O^t_{jjmcps}$ is the tonnage of commodity $c$ shipped from county $j \in J^{\prime}$ to county $m \in J^{{\prime}{\prime}}$ via truck in month $p$ in scenario $s$
    \item $O^r_{jmcps}$ is the tonnage of commodity $c$ shipped from county $j \in J^{\prime}$ to county $m \in J^{{\prime}{\prime}}$ via rail in month $p$ in scenario $s$
    \item $Q_{jcps}$ is the tonnage of commodity $c$ shortage in county $j \in J^{{\prime}{\prime}}$ in month $p$ in scenario $s$
\end{itemize}
\par
\noindent Our proposed mathematical model is defined as follows.
 \begin{subequations}
 \small
 \begin{align}\label{equ:firststage}
    (WSN): min\hspace{0.10 in}\sum_{i \in I, e \in E}\kappa_{e}Z_{ie} + \sum_{i \in I, f \in F}\iota_fY_{if} + E(H(X,\tilde{d}))
    \end{align}
Subject to:
        \begin{align}
        \sum_{i\in I}\sum_{e\in E}\kappa_{e}Z_{ie} +\sum_{i\in I}\sum_{f\in F} \iota_fY_{if} \leq b\label{const:budget}\\
        	Z_{ie} \in Z^+,\hspace{0.5in}\forall i \in I, e \in E,\label{cons:Zinteger}\\
	 Y_{if} \in Z^+,\hspace{0.5in}\forall i \in I, f \in F. \label{cons:Yinteger}
	\end{align}
 \end{subequations}
The function (\ref{equ:firststage}) represents the objective function that consists of the total of first stage costs (capacity expansion) and the expected second stage cost (transportation costs). The first term of (\ref{equ:firststage}) represents the total cost of new operational equipment (loading / unloading) and the second term represents the total cost of new storage facilities. Constraint (\ref{const:budget}) ensures that the total cost does not exceed the total available budget allocated for investments to improve the inland waterway system. Constraints (\ref{cons:Zinteger}) and (\ref{cons:Yinteger}) determine the integer variables related to the decisions on storage facilities and operational equipment, respectively. 
Let $\mathcal{Y}=\{Z_{ie}, Y_{if}| i\in I, e\in E, f\in F\}$ represent solutions to \eqref{equ:firststage} to \eqref{cons:Yinteger}. For a given value $\mathcal{\bar{Y}}\in \mathcal{Y}$, and a realization $d$ of random demand $\tilde{d}$, the following is the formulation of the second stage problem $H(\mathcal{\bar{Y}},d)$.\par

\begin{subequations}
\small
\begin{align}
\begin{split}
\label{equ:secondstagecost}H(\mathcal{\bar{Y}},d) & = min \sum_{c\in C } \sum_{p\in P}
  \left(\sum_{i\in I^{\prime}}\sum_{j\in J^{\prime}}\sum_{e \in E}( \alpha^t_{ij}X^t_{jicpes} + \alpha ^r_{ij}X^r_{jicpes}) + 
\sum_{i \in I^{\prime}}\sum_{f \in F}h_{c}U_{icpf}+\sum_{i\in I^\prime}\sum_{k \in I^{\prime \prime}}a_{ik}W_{ikcpe}+\right.\\
&\hspace{0.18 in}\left. \sum_{i \in I^{{\prime}{\prime}}}\sum_{f \in F}h_{c}V_{icpf}+  \sum_{i \in I^{{\prime}{\prime}}}\sum_{j \in J^{{\prime}{\prime}}}\sum_{e \in E}( \alpha^t_{ij}R^t_{ijcpe}+ \alpha ^r_{ij}R^r_{ijcpe})+ \right.\\
&\hspace{0.18 in}\left. \sum_{j \in J^\prime}\sum_{m \in J^{\prime \prime}}(l^t_{jm}O^t_{jmcp}+
l^r_{jm}O^r_{jmcp})+
 \sum_{j \in J^{{\prime}{\prime}}}\mu*Q_{jcp}\right)
\end{split}
\end{align}

\noindent subject to:

\small
\begin{align}
    \label{const:supply}\sum_{m \in J^{\prime \prime}}(O^t_{jmcp}+O^r_{jmcp})+\sum_{i\in I^{\prime}}\sum_{e\in E}(X^t_{jicpe}+X^r_{jicpe}) \leq q_{jcp} \hspace{0.3in} \forall j \in J^{{\prime}}, c \in C, p \in P
\end{align}
\begin{align}
    \label{cons:shortage}Q_{mcp}+\sum_{i\in I^{{\prime}{\prime}}}\sum_{e\in E}(R^t_{imcpe}+R^r_{imcpe}) + \sum_{j \in J^{\prime}}(O^t_{jmcp}+O^r_{jmcp}) = d_{mcp}
    \hspace{0.3in}  \forall m \in J^{{\prime}{\prime}}, c \in C, p \in P
\end{align}
\begin{align}
    \label{cons:FBO}\sum_{j\in J^{\prime}}\sum_{e\in E}(X^t_{jicpe}+X^r_{jicpe}) + \sum_{f \in F}U_{icp-1f}= \sum_{f \in F}U_{icpf} + \sum_{k \in I^{{\prime}{\prime}}}\sum_{e\in E}W_{ikcpe} 
    \hspace{0.3in} \forall i \in I^{\prime}, c \in C, p \in P
\end{align}
\begin{align}
    \label{cons:FBD}\sum_{i\in I^{\prime}}\sum_{e\in E} W_{ikcpe} + \sum_{f \in F}V_{kcp-1f}= \sum_{f \in F}V_{kcpf} + \sum_{m \in J^{{\prime}{\prime}}}\sum_{e\in E}(R^t_{km{cpe}}+R^r_{kmcpe}) 
    \hspace{0.3in} \forall k \in I^{{\prime}{\prime}}, c \in C, p \in P
\end{align}
\begin{align}
    \label{cons:processing}\sum_{c\in C}\left(\Lambda_{ec}\sum_{j\in J^{{\prime}{\prime}}}(R^t_{ijcpe}+R^r_{ijcpe}) +  \sum_{k \in I^{{\prime}{\prime}}}(W_{ikcpe}+W_{kicpe})+ \nonumber\right.
    \left.\sum_{j\in J^{\prime}} (X^t_{jicpe}+X^r_{jicpe})\right) \leq \ m_e(n_{ie}+\bar{Z}_{ie})  \\
    \forall i \in I, p \in P,e \in E
\end{align}
\begin{align}
\label{cons:inventory}\sum_{c\in C}\zeta_{fc}(U_{icpf} + V_{icpf}) \leq l_{f}(k_{if}+\bar{Y}_{if})\hspace{0.3in} \forall i \in I, p \in P,f \in F      
\end{align}
\begin{align}
   \label{cons:UInitial}U_{icpf}=0 \hspace{0.6in} \forall i \in I^{\prime}, c \in C, p \in \{0\}, f \in F 
\end{align}
\begin{align}
    \label{cons:VInitial}V_{icpf}=0 \hspace{0.6in} \forall i \in I^{{\prime}{\prime}}, c \in C, p \in \{0\}, f \in F 
\end{align}
\begin{align}
    \label{cons:nonnegativity}X^t_{jicpe},X^r_{jicpe}, U_{icpf}, W_{ikcpe} ,V_{icpf}, R^t_{ijcpe}, R^r_{ijcpe}, O^t_{jmcp},
    O^r_{jmcp}, 
    Q_{jcp} \in R^+ \hspace{0.6in} \forall i \in I^\prime, k \in I^{{\prime}{\prime}}, j \in J^\prime, \nonumber \\
     m \in J^{{\prime}{\prime}}, c \in C, p \in P, f \in F
\end{align}
\end{subequations}
Function (\ref{equ:secondstagecost})  minimizes the total supply chain costs. These costs include shipping cost via truck, rail, and barges, inventory cost, and a penalty cost for any unmet demand. Constraints (\ref{const:supply}) ensure that the amount of commodity shipped does not exceed the quantity of supply available.  Constraints (\ref{cons:shortage}) capture the shortage of commodity in cases when demand exceeds supply. Constraints (\ref{cons:FBO}) and (\ref{cons:FBD}) are the flow balance constraints for origin and destination ports. Constraints (\ref{cons:processing}) ensure that the commodity handled at the origin and destination ports does not exceed the port capacity. Constraints (\ref{cons:inventory}) ensure that the total inventory of a port does not exceed the inventor holding capacity of that port. Constraints (\ref{cons:UInitial}) and (\ref{cons:VInitial}) assign initial inventory at the destination and origin ports, respectively as zero. Constraints (\ref{cons:nonnegativity}) are the non-negativity constraints.   

\section{Solution Approach}
The computational burden of solving the problem formulated in (\ref{equ:firststage})-(\ref{cons:nonnegativity}) warrants the use of various solution approaches. This section discusses the Benders Decomposition algorithm and techniques to accelerate its convergence rate.  

\subsection{Benders Decomposition Algorithm}

\noindent The uncertainty in commodity demand (Figure \ref{fig:stochasticdemand}) requires examining a number of scenarios for a robust network design. We approximate the distribution of our stochastic demand via a discrete distribution. Let $S$ represent the discrete set of demand realization and let $\omega_{s}$ $\forall s \in S$ represent the corresponding probabilities.\par
The solution to the 2-SOP problem can be computationally expensive based on the size of the problem given by $|I|$, $|J|$,$|C|$,$|P|$ and $|S|$. To overcome this computational burden, Benders decomposition algorithm \cite{benders1962partitioning}, widely used to solve large-size, mixed integer linear problems, is employed. Following this method, we decompose original problem into two subproblems: an integer master problem (MP) (\ref{equ:MWSN}) and $|S|$ linear subproblems (\ref{equ:secondstage}). The MP along with an auxiliary variable and optimality cut provides an approximation of the original problem. Benders decomposition is an iterative procedure. Let $\mathcal{\bar{Y}}$ be the solution of MP (\ref{equ:MWSN}). For the given $\mathcal{\bar{Y}}$, $|S|$ subproblems are solved, one for each realization $d_{s}$ of the stochastic demand $\tilde{d}$. Solutions to the subproblems are used to develop feasibility and optimality cuts that are added to the MP. These cuts ensure that, if the current solution $\mathcal{\bar{Y}}$ of the MP is not feasible or optimal to the original \textit{(WSN)} problem, this solution is excluded from the feasible region and will not be used in other iterations of the Benders algorithm. In each iteration of the algorithm, a lower bound and an upper bound are generated. Objective value of $(M-WSN)$ provides a lower bound, and the solution of $(M-WSN)$ and $(S-WSN(s))$ is used to calculate an upper bound for the original problem $(WSN)$ This is continued until the relative gap between the lower and upper bound converge to a given threshold. The following model (\ref{equ:MWSN}) is the MP.

\begin{subequations}\label{equ:MWSN}
\small
\begin{align}
    {(M-WSN): min \sum_{i \in I e \in E}\kappa_{e}Z_{ie} + \sum_{i \in I f \in F}\iota_fY_{if} +\sum_{s \in S}\omega_{s} \theta_{s}}
\end{align}

\noindent subject to: (\ref{const:budget})-(\ref{const:intY})

\small
\begin{align}
\begin{split}
\label{cons:theta}
\theta^{n}_{s} \geq \sum_{j \in J} \sum_{c \in C} \sum_{p \in P} (\nu_{jcps}q_{jcp}+\xi_{jcps}d_{jcps}+
\sum_{i \in I}\sum_{e \in E}\sum_{p \in P}\pi_{ipes}m_{e}(n_{ie}+Z_{if})+\sum_{i \in I}\sum_{f \in F}\sum_{p \in P}\sigma_{ipfs}l_{f}(k_{if}+Y_{if})  \\ \forall s \in S, \hspace{0.1in} n = \{1, \ldots, N^{\prime}\}
\end{split}
\end{align}
\end{subequations}
where $N^{\prime}$ is the current number of iterations.

\noindent The Benders decomposition solves $(M-WSN)$ iteratively and in each iteration $n$, let $\mathcal{\bar{Y}}^{n}$ represent the corresponding solution. Let $\Phi^{n}$ be the objective function value of model (\ref{equ:MWSN}) obtained at the $n^{th}$ iteration. Let $\phi^{n}$, given by equation (\ref{equ:z_IF}),  be the corresponding cost of infrastructure investment. 

\begingroup
\small
\begin{align}
    \label{equ:z_IF}\phi^{n} = \sum_{i \in I e \in E}\kappa_{e}Z^{n}_{ie} + \sum_{i \in I f \in F}\iota_fY^{n}_{if}
\end{align}
\endgroup

\noindent Given $\mathcal{\bar{Y}}^{n}$, following scenario-based subproblem ($S-WSN(s)$) is solved for each scenario $s \in S$.

\small
\begin{align}
\label{equ:secondstage}(S-WSN(s)) :& min \sum_{c\in C} \sum_{p\in P}
  (\sum_{i\in I^{\prime}}\sum_{j\in J^{\prime}}\sum_{e\in E}( \alpha^t_{ij}X^t_{jicpes} + \alpha ^r_{ij}X^r_{jicpes}) + 
\sum_{i \in I^{\prime}}h_{c}U_{icpfs}+\sum_{i\in I^{\prime}}\sum_{k \in I^{{\prime}{\prime}}}\sum_{e\in E}a_{ik}W_{ikcpes}+\\\nonumber
&\hspace{0.18 in}\sum_{i \in I^{{\prime}{\prime}}}h_{c}V_{icpfs}+  
\sum_{i \in I^{{\prime}{\prime}}}\sum_{j \in J^{{\prime}{\prime}}}\sum_{e\in E} ( \alpha^t_{ij}R^t_{ijcpes}+ \alpha ^r_{ij}R^r_{ijcpes})+
\sum_{j \in J^{\prime}}\sum_{m\in J^{{\prime}{\prime}}}(l^t_{jm}O^t_{jjmcps}+
l^r_{jm}O^r_{jmcps})+\\
&\hspace{0.18 in}\sum_{j \in J^{{\prime}{\prime}}}\mu*Q_{jcps})\nonumber
\end{align}
\hspace{1in} Subject to:  (\ref{const:supply})-(\ref{cons:nonnegativity})\newline 

The constraint (\ref{cons:theta}) is the scenario specific optimality cut added to ($M-WSN$) in each iteration of the Benders decomposition algorithm. $\nu_{jcps} \forall j \in J^\prime$,  $\xi_{jcps} \forall j \in J^{\prime \prime}$, $\upsilon_{icps}$, $\chi_{icps}\forall i \in I^{\prime \prime}$, $\pi_{ipes}$, $\sigma_{ipfs}$ are dual variables for constraints (\ref{const:supply}) to (\ref{cons:inventory}) respectively. Let $\mathcal{X}^n(s)$ denote these solutions for scenario $s$ in $n^{th}$ iteration. Constraint (\ref{cons:shortage}) makes the scenario-based subproblems $(S-WSN(s))$ always feasible for any values in first-stage decision variables. This constraint ensures that the demand is satisfied either by the supply from counties within the state or via other states. For this reason, we do not need to add feasibility cuts to the MP. The dual problem of $(S-WSN(s))$ is shown in Appendix \ref{appendix:dual}.

Let $\Theta^{n}_s$ be the objective value of the function $s^{th}$ $(S-WSN(s))$ for scenario $s$. We calculate $\Theta^{n}$ as follows. 

\begingroup
\small
\begin{align}
    \label{equ:z_SP}\Theta^{n} = \sum_{s \in S}\omega_{s}\Theta^{n}_{(s)}
\end{align}
\endgroup
The pseudo-code of the Benders Decomposition algorithm is shown in Algorithm (\ref{alg:Benders}).

\begin{algorithm}
\caption{Benders Decomposition Algorithm}\label{alg:Benders}
Initialize $\epsilon$. Set $n \gets 1$, $LB^{n} \gets -\infty$, $UB^{n} \gets +\infty$, $abort \gets false$\\
\While{abort = false}{
  Solve $(M-WSN)$ to obtain $\Phi^{n}$, $\phi^{n}$ and $\mathcal{Y}^{n}$\\
  \If{$\Phi^{n}$ > LB}{
    $LB^{n} \gets \Phi^{n}$\\
  }
  For all $s \in S$, Solve $(S-WSN(s))$ to obtain $\Theta^{n}_{s}$ and $\mathcal{X}^{n}(s)$\\
  \If{$UB^{n} > \Theta^{n} +\phi^{n}$}{
      $UB^{n} \gets \Theta^{n} +\phi^{n}$\\
    }
  \eIf{$\frac{UB^{n}-LB^{n}}{UB^{n}} \leq \epsilon$}{
  $abort \gets true$}{
  Add cut (\ref{cons:theta}) to $(M-WSN)$\\
  $n \gets n+1$\\
  }
}
return $UB$, $\mathcal{Y}^{n}$ and $\mathcal{X}^{n}(s)$
\end{algorithm}

\subsection{Methods to Accelerate Benders Decomposition Algorithm}
The Benders decomposition algorithm is known to be extremely slow and computationally expensive. \citep{mcdaniel1977modified}, \citep{cote1984large}, \citep{poojari2009improving}, \citep{magnanti1981accelerating}, \citep{rei2009accelerating}, \citep{saharidis2010accelerating} and \citep{saharidis2011initialization} have developed and successfully implemented acceleration techniques to enhance the Benders Decomposition algorithm. In our study we add Knapsack inequalities and Pareto-optimal cuts to improve the convergence rate of the algorithm. The details of the proposed techniques are provided in later parts of this section.\par
\noindent \emph{Knapsack inequalities:} Adding knapsack inequalities 
can improve the convergence rate of the algorithm by reducing  the solution space of (M-WSN), thus, reducing the time it takes to solve the  problem \cite{santoso2005stochastic}.  In addition, commercial solvers such as GUROBI and CPLEX, can extract valid inequalities from knapsack inequalities. 

Let $LB^{n}$ denote the lower bound obtained in the $n^{th}$ iteration of the algorithm. Therefore, the following knapsack inequalities are added to ($M-WSN$) in the $(n+1)^{th}$ iteration to accelerate the convergence rate of the Benders decomposition algorithm: 
\begin{eqnarray}
\begin{split}
\label{cons:knapackLB}LB^{n}\leq \sum_{i \in I e \in E}\kappa_{e}Z_{ie} + \sum_{i \in I, f \in F}\iota_fY_{if} +\sum_{s \in S}\omega_{s}\theta_{s}.
\end{split}
\end{eqnarray}

\par
\noindent \emph{Pareto-optimal Cuts:} The subproblems ($S-WSN(s)$) are capacitated transportation problems. The transportation problem is degenerate in nature \cite{ahuja1988network}, that is, it has multiple optimal solutions and each solution generates optimality cuts of different strength.  Hence, the solution of the subproblems should be chosen in such a way that it produces the strongest cuts. \citep{magnanti1981accelerating} found that adding Pareto-optimal cuts to the MP improves the convergence rate of the Benders Decomposition algorithm. The generation of Pareto-optimal cuts proposed by \citep{magnanti1981accelerating} requires solving two subproblems, one associated with solutions of the MP and another associated with core points. A point $y \in ri(Y^{c})$ is called a core point of $Y$, where $ri(S)$ and $S^{c}$ are the relative interior and convex hull of set $S \subseteq \mathbb{R}^{k}$ respectively \cite{papadakos2008practical}. Let, $\bar{\mathcal{Y}}^{n}$ be the solution of the MP at $n^{th}$ iteration and 
$\mathcal{\bar{Y}}^{o} = \{Z^{o,n}_{ie}=0, Y^{o,n}_{if}=0|n=1,i \in I, e \in E,  f \in F\}$ 
be the set of initial core points. The subproblems to generate Pareto-optimal cuts are given as:

\begin{subequations}
\begin{align}
\begin{split}
\label{equ:MWdualsecondstage}(MWS-WSN(s)): & \max \sum_{j \in J}\sum_{c\in C}\sum_{p \in P}\sum_{s \in S}(\nu_{jcps}q_{jcps} + \xi_{jcps}d_{jcps}) +
\sum_{i \in I}\sum_{e\in E}\sum_{p \in P}\sum_{s \in S}\pi_{ipes}m_{e}(n_{ie}+Z^{o}_{ie})+\\
&\hspace{0.18 in}\sum_{i \in I}\sum_{f\in F}\sum_{p \in P}\sum_{s \in S}\sigma_{ipfs}l_f(k_{if}+Y^{o}_{if})
\end{split}
\end{align}

\noindent subject to (\ref{cons:DS_xt})-(\ref{cons:DS_free})

\small
\begin{align}
\begin{split}
\label{cons:MWSP}\sum_{j \in J}\sum_{c\in C}\sum_{p \in P} \left(\nu_{jcps}q_{jcps} + \xi_{jcps}d_{jcps}\right) + 
\sum_{i \in I}\sum_{e\in E}\sum_{p \in P}\pi_{ipes}m_{e}(n_{ie}+Z^{n}_{ie})+\\
\sum_{i \in I}\sum_{f\in F}\sum_{p \in P}\sigma_{ipfs}l_f(k_{if}+Y^{n}_{if}) = \Theta^{n}(s)\\ 
\hspace{0.2in} \forall s \in S
\end{split}
\end{align}
\end{subequations}

Since this technique relies on the solution of the subproblems, \citep{papadakos2008practical} proposed a methodology to generate sub problem independent Pareto-optimal cuts. \citep{papadakos2008practical} showed that by using different core points in each iteration, the constraint (\ref{cons:MWSP}) could be ignored. Pareto-optimal cuts generated by this method are known as \textit{modified Magnanti-Wong pareto-optimal cuts}. The core points are updated in every iteration as follows.

\begingroup
\small
\begin{align}
\label{equ:UpdateZCore}Z^{o,n+1}_{ie} = (1-\lambda)Z^{o,n}_{ie}+\lambda Z^{n}_{ie}\\
\label{equ:UpdateYCore}Y^{o,n+1}_{if} = (1-\lambda)Y^{o,n}_{if}+\lambda Y^{n}_{if}
\end{align}
\endgroup
\citep{papadakos2008practical} and \citep{mercier2005computational} empirically showed that the value of $\lambda$ = 0.5 gives the best result. Pseudo code of the accelerated Benders Decomposition with Pareto-optimal cuts is provided in Algorithm (\ref{alg:BenderswithParetoOptimal}).
\RestyleAlgo{ruled}

\begin{algorithm}[tbh]
Initialize, $\epsilon$, $Z^{o,n}_{ie}$,$Y^{o,n}_{if}$,Set, $n \gets 1$,$LB^{n} \gets -\infty$,$UB^{n} \gets +\infty$,$abort \gets false$\\
\While{$abort = false$}{
  For all $s \in S$, Solve $(MWS-WSN(s))$ to obtain $\Theta^{n}_{s}$ and $\mathcal{Y}^{n}(s)$ \\
  Add cuts (\ref{cons:theta}) to $(M-WSN)$\\
  Solve $(M-WSN)$ to obtain $\Phi^{n}$, $\phi^{n}$ and $\mathcal{Y}^{n}$\\
  \If{$\Phi^{n} > LB$}{
    $LB^{n} \gets \Phi^{n}$\\
  }
  For all $s \in S$, Solve subproblem $(S-WSN(s))$ to obtain $\Theta^{n}_{s}$ and $\mathcal{X}^{n}(s)$\\
  \If{$UB^{n} > \Theta^{n} + \phi^{n}$}{
      $UB^{n} \gets \Theta^{n} + \phi^{n}$\\
    }
  \eIf{$\frac{UB^{n}-LB^{n}}{UB^{n}} \leq \epsilon$}{
  $abort \gets true$}{
  Add cuts (\ref{cons:theta}) to $(M-WSN)$\\
  Add Knapsack inequalities (\ref{cons:knapackLB}) to $(M-WSN)$\\
  Update core points using equations (\ref{equ:UpdateZCore}) and (\ref{equ:UpdateYCore})\\
  $n \gets n+1$ 
  }
  }
\KwRet{$UB, \mathcal{Y}^{n}, \mathcal{X}^{n}(s)$}
\caption{Benders Decomposition Algorithm with Knapsack inequalities and Pareto-optimal Cuts}\label{alg:BenderswithParetoOptimal}
\end{algorithm}

\section{Computational study and managerial insights}
This section discusses the solution quality of the proposed algorithm, a case study application of the model to the Arkansas Section of the MKARNS, and the performance evaluation of stochastic solutions.  
\subsection{Performance evaluation}
We solved a variety of problems to determine the quality of the proposed algorithms (Table \ref{tab:Case Size}) and compare five solution approaches: $(1)$ Gurobi solver, $(2)$ Benders algorithm, $(3)$ Benders with Knapsack inequalities, $(4)$ Benders with Knapsack inequalities and Pareto-optimal cuts. To terminate the algorithm, we consider following stopping criterion: $(i)$ optimality gap $\leq 1\%$, $(ii)$ number of iterations $\geq 500$ and $(iii)$ algorithm run time $\geq$ 12,600 seconds. The experiments are carried out on a Windows 10 PC with an Intel Core i7 3.2 GHz processor and 32 GB of RAM. The results of the experiments are summarized by the algorithm run time in seconds $t(s)$, the optimality gap $\epsilon(\%)$, and the number of iterations (n) at the time when the stop criteria is met (Table \ref{tab:sol_comparision}). \par
\begin{table}
\caption{Test Case Size} \label{tab:Case Size}
\centering
\renewcommand{\arraystretch}{1.2}
\begin{tabular}{c | c c c c c | c |c c c c c}
    \hline
    Problem Nr. & $|I|$& $|J|$& $|C|$& $|P|$& $|S|$& Problem Nr. & $|I|$& $|J|$& $|C|$& $|P|$& $|S|$\\
    \hline
1	&	15	&	45	&	5	&	6	&	4	&	9	&	30	&	75	&	11	&	9	&	12	\\
2	&	15	&	45	&	5	&	6	&	8	&	10	&	30	&	75	&	11	&	9	&	16	\\
3	&	15	&	45	&	5	&	6	&	10	&	11	&	30	&	75	&	11	&	12	&	4	\\
4	&	15	&	45	&	5	&	6	&	12	&	12	&	30	&	75	&	11	&	12	&	8	\\
5	&	15	&	45	&	5	&	6	&	16	&	13	&	30	&	75	&	11	&	12	&	10	\\
6	&	30	&	75	&	11	&	9	&	4	&	14	&	30	&	75	&	11	&	12	&	12	\\
7	&	30	&	75	&	11	&	9	&	8	&	15	&	30	&	75	&	11	&	12	&	16	\\
8	&	30	&	75	&	11	&	9	&	10	&		&		&		&		&		&		\\
    \hline
    
\end{tabular}

\end{table}

Overall, Gurobi outperforms the alternative solution approaches in terms of run time, optimality gap, and number of iterations in all experiments except for the cases (8)-(10) and (12)-(15). This warrants the implementation of the Benders decomposition algorithm. Although the standard Benders decomposition does not have memory issues, it fails to converge to solutions with a smaller than 1\% optimality gap within the given time limit for medium and large sized problems (problems (6) - (15), Table \ref{tab:sol_comparision}). This finding warrants the application of acceleration techniques to enhance the algorithm. With the addition of Knapsack inequalities, improvement in terms of solution time, optimality gap, and number of iterations is observed for most of the cases. 
However, similar to the standard Benders decomposition, it is not successful in finding solutions with smaller than 1\% optimality gap within 12,600 seconds (3.5 hours). 
We observe a significant improvement in terms of solution quality and run time  when we add both Knapsack inequalities and Pareto-optimal cuts to the MP. Benders decomposition algorithm with Knapsack inequalities and Pareto-optimal cuts provides solutions of high quality within the proposed threshold run time. 

\begin{table}[!htb]
\small
\caption{Performance evaluation of solution approaches}
\label{tab:sol_comparision}
  \renewcommand{\arraystretch}{1.2}
    \begin{tabular}{ccccccccccccccc}
    \hline
    \multirow{2}{*}{Case} & 
    \multicolumn{2}{c}{Gurobi} &       
    & \multicolumn{3}{c}{Benders} &       
    & \multicolumn{3}{c}{Benders + KI} &      
    & \multicolumn{3}{c}{Benders + KI + PO} \\
    \cline{2-3}
    \cline{5-7}
    \cline{9-11}
    \cline{13-15}
          & t (s) & $\epsilon$ (\%) &       & t (s) & $\epsilon$ (\%) & n     &       & t (s) & $\epsilon$ (\%) & n     &       & t (s) & $\epsilon$ (\%) & n \\
    \hline
    1     & \textbf{             5 } & 0.90  &       & 40    & 0.51  & 21    &       & 38    & 0.45  & 19    &       & 26    & 0.95  & 7 \\
    2     & \textbf{           20 } & 0.07  &       & 160   & 0.99  & 38    &       & 135   & 0.98  & 31    &       & 119   & 0.28  & 15 \\
    3     & \textbf{           25 } & 0.16  &       & 378   & 0.99  & 64    &       & 349   & 1.00  & 55    &       & 131   & 0.85  & 13 \\
    4     & \textbf{           36 } & 0.22  &       & 185   & 0.97  & 29    &       & 151   & 0.98  & 24    &       & 144   & 0.92  & 12 \\
    5     & \textbf{           52 } & 0.09  &       & 161   & 0.98  & 20    &       & 161   & 0.97  & 20    &       & 157   & 0.99  & 10 \\
    6     & \textbf{         955 } & 0.51  &       &        12,600  & 5.09  & 179   &       &        12,600  & 4.78  & 135   &       &          2,960  & 0.81  & 30 \\
    7     & \textbf{      4,063 } & 0.46  &       &        12,600  & 4.96  & 118   &       &        12,600  & 4.51  & 99    &       &          5,784  & 0.88  & 30 \\
    8     &  -    & -     &       &        12,600  & 4.52  & 95    &       &        12,600  & 4.41  & 88    &       & \textbf{         5,442 } & 0.96  & 23 \\
    9     &  -    & -     &       &        12,600  & 3.24  & 87    &       &        12,600  & 2.76  & 80    &       & \textbf{         6,702 } & 0.48  & 24 \\
    10    &  -    & -     &       &        12,600  & 2.15  & 66    &       &        12,600  & 2.02  & 63    &       & \textbf{         7,421 } & 0.89  & 20 \\
    11    & \textbf{      1,732 } & 0.77  &       &        12,600  & 5.48  & 157   &       &        12,600  & 5.20  & 135   &       &          5,232  & 0.80  & 36 \\
    12    &  -    & -     &       &        12,600  & 5.33  & 88    &       &        12,600  & 5.10  & 83    &       & \textbf{         8,820 } & 0.98  & 31 \\
    13*    &  -    & -     &       &        12,600  & 5.41  & 72    &       &        12,600  & 4.94  & 67    &       & \textbf{       10,245 } & 0.94  & 29 \\
    14    &  -    & -     &       &        12,600  & 3.77  & 61    &       &        12,600  & 3.42  & 58    &       & \textbf{       10,407 } & 0.85  & 25 \\
    15    &  -    & -     &       &        12,600  & 2.70  & 46    &       &        12,600  & 2.70  & 46    &       & \textbf{       12,586 } & 0.69  & 23 \\

    \hline
    \multicolumn{15}{l}{\small - Out of memory} \\
    \multicolumn{15}{l}{\small * Representative case size for a real-world problem}\\
    \end{tabular}%

\end{table}

\subsection{Case Study: McClellan-Kerr Arkansas River Navigation System}
\begin{figure}[h!tb]
  \centering
  \includegraphics[width=0.8\textwidth]{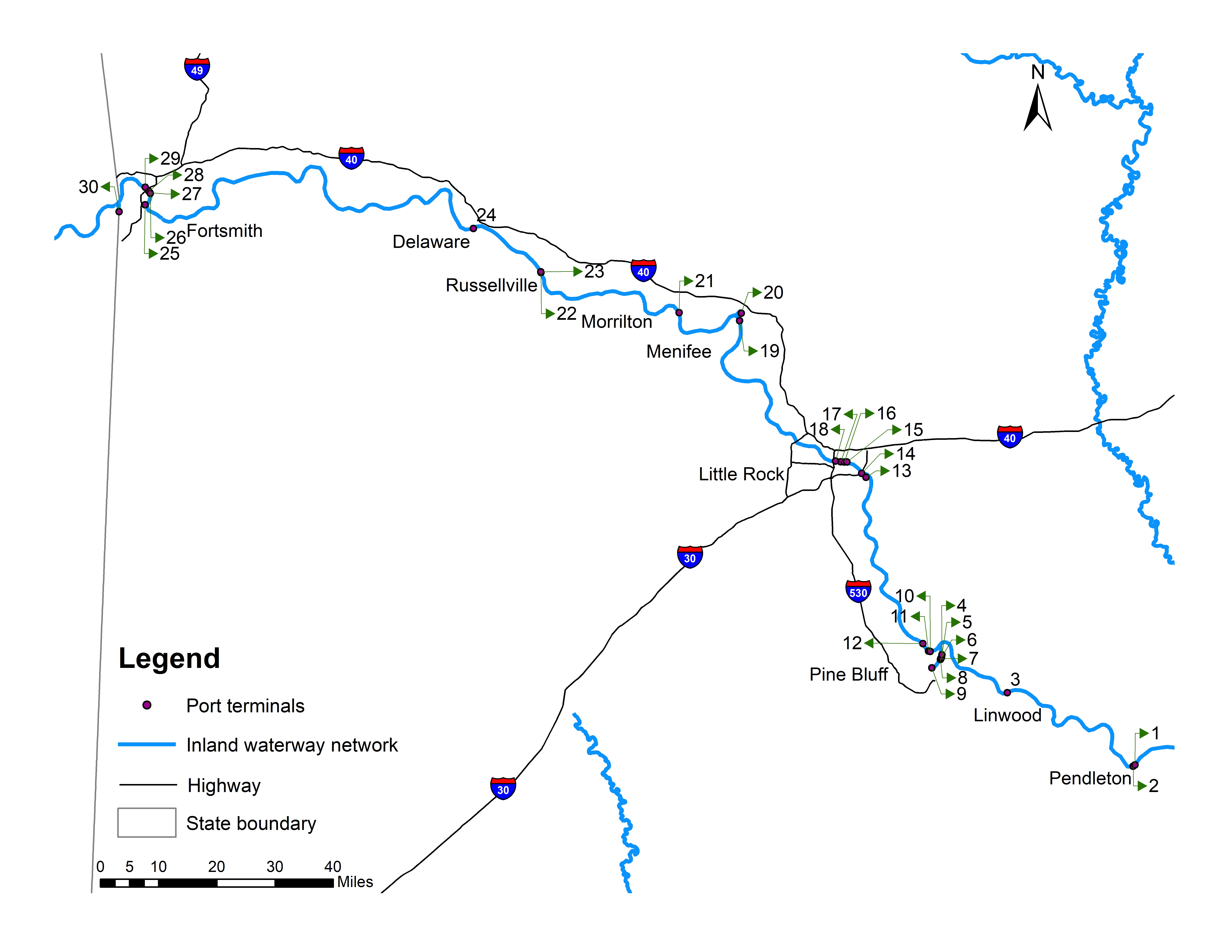}
  \caption{Arkansas Segment of MKARNS.}\label{fig:waterwaymap}
\end{figure}

The 2-SOP problem is applied to the Arkansas section of the MKARNS. This 308 mile inland waterway system has 13 locks and 43 freight ports (Figure \ref{fig:waterwaymap}). The Arkansas segment of MKARNS has significant economic impact in Arkansas in terms of sales, Gross Domestic Product, labor income, and jobs \cite{nachtmann2015regional}.\par  
\subsubsection{Data Description}
\emph{Commodity demand and supply data}\\
The US Army Corps of Engineers (USACE) produces the Lock Performance Management System (LPMS) report which records monthly commodity volume (in tons) passing
through each lock. The LPMS classifies commodities into 9 groups and several subgroups. Data is available for the years 2009 through 2016. This source was used to create demand scenarios in the 2-SOP model. \par
The Arkansas Statewide Travel Demand Model \cite{STDM2012} maintained by the Arkansas Department of Transportation (ARDOT) reports amounts (in tons) of commodities shipped between counties. The total amount of commodity shipped from a county is considered to be the supply for the (shipment) origin county. Similarly, the total amount of commodities received by a county is considered to be the demand for that (shipment) destination county. The commodity grouping in the Statewide Travel Demand Model and LMPS differed.  We reorganized commodities into 11 groups that align with the subgroups reported by LMPS.\par
\noindent\emph{Transportation cost data}\\
The transportation cost (in \$/ton-mile) for truck and rail are derived from the ARDOT Travel Demand Model \cite{STDM2012}, \citep{BTSTransportationStatistics2010} and \citep{SurfaceBoard2003} . The barge transportation cost (in \$/ton-mile) is derived from the data reported by \citep{GrainTransportationReport} for major cities along the waterways. However, the data is not reported for cities along the Arkansas River. Thus, rates from the surrounding region in St. Louis and Cairo-Memphis were averaged to estimate the costs for the Arkansas river. Transportation cost data used in our study are shown in Table \ref{tab:transportationcost} in Appendix \ref{appendix:a}. \par
\noindent\emph{Port capacity data}\\
The National Transportation Atlas Database (NTAD), published by the BTS\cite{BTSPortsShapefile} contains information on the location, commodity handling equipment, storage facility and road / rail connection of terminals at US coastal, Great Lakes, and inland ports. Data about the number of equipment units and storage facilities  are gathered from this dataset. These data are supplemented and updated using information obtained from ports' websites and satellite images. The data for port processing and storage capacity used in our study are shown in Table \ref{tab:processingcapacity} and \ref{tab:storagecapacity}, respectively, in Appendix \ref{appendix:a}.\par
\noindent\emph{Infrastructure costs}\\
Costs for equipment (i.e., crane, conveyor, hooper, forklift) and storage facilities (i.e., warehouse, storage tank, paved and unpaved storage) are obtained from \citep{braham2017locating}. These costs include labor and materials, as well as general overhead. \citep{braham2017locating} selected these costs from a material, construction, and equipment cost database from \citet{RSMeans2014,RSMeans2017} and validated through interviews with industry representatives.  All costs for our study are calculated based on the 2020 dollar value. The infrastructure costs used in our study are shown in Table \ref{tab:equipmentcost} and \ref{tab:storagefacilitycost} in Appendix \ref{appendix:a}. \par
\noindent\emph{Scenario definitions}\\
To capture the impacts of stochasticity in commodity demand to infrastructure decisions, we generate 10 different demand scenarios. Scenarios 1 to 8 are based on historical commodity throughput data gathered from the LPMS between 2009 and 2016. We assign a probability of 6.5\% to the scenario developed using historical data from 2009, and increase this probability by 0.5\% per year for subsequent years. This is done with the assumption that the demand scenarios developed from recent years have more probability of occurrence compared to scenarios developed from previous years. The probability of scenarios developed from historical data sum to 66\%. Two additional scenarios: scenario 9 with probability of 20\% and scenario 10 with probability of 14\%, are developed based on the 15 and 25 year commodity demand projection from BTS \cite{FreigthFacts2019}, respectively.

\subsubsection{Case study results}
\begin{figure}[!htb]
    \centering
    \begin{minipage}[t]{.5\textwidth}
      \begin{subfigure}{\linewidth}
      \centering
      \includegraphics[width=1\linewidth]{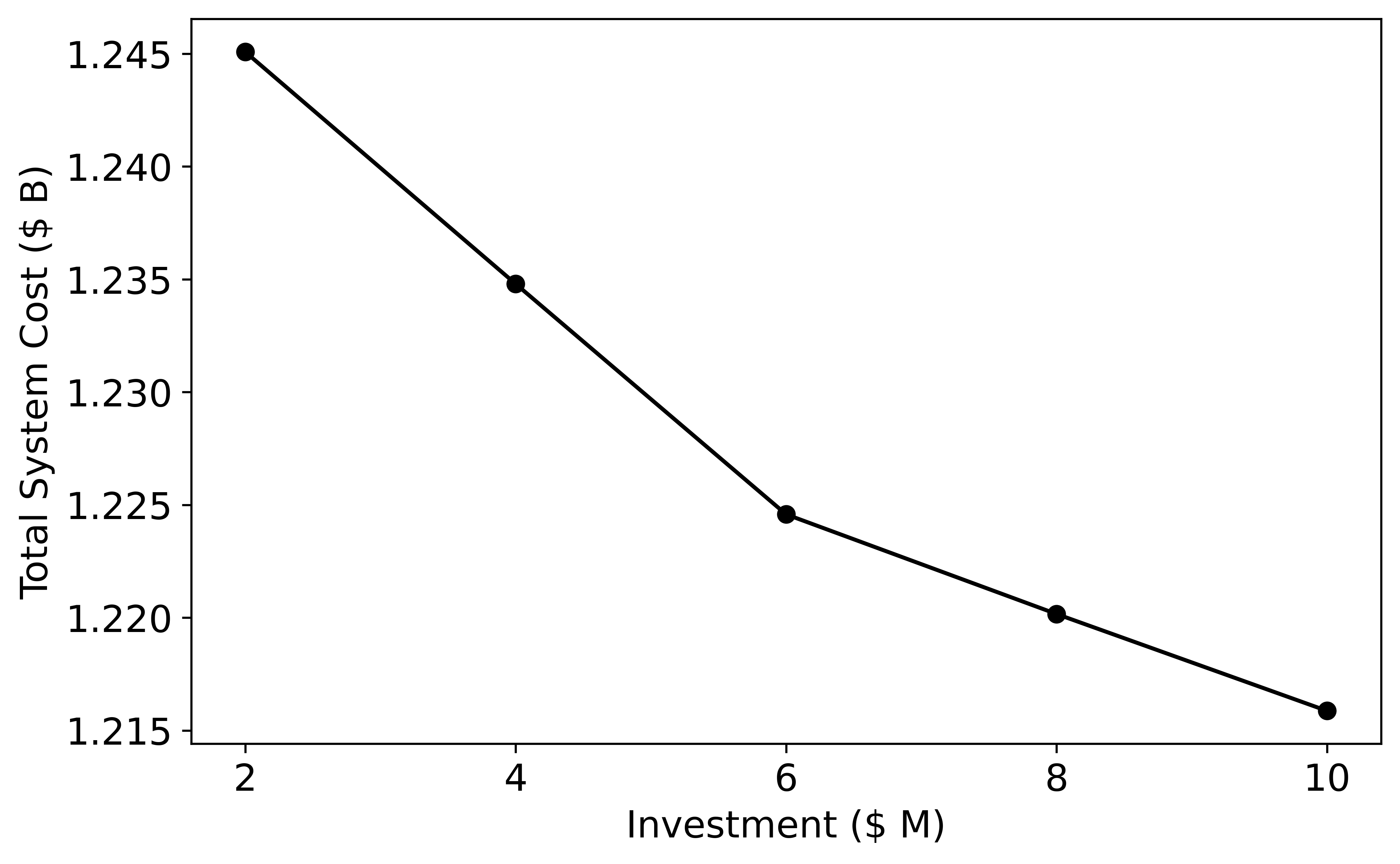}
      \caption{System Cost}
      \label{fig:TotalSystemCost}
    \end{subfigure}
    \end{minipage}%
    \begin{minipage}[t]{.5\textwidth}
      \begin{subfigure}{\linewidth}
      \centering
      \includegraphics[width=1\linewidth]{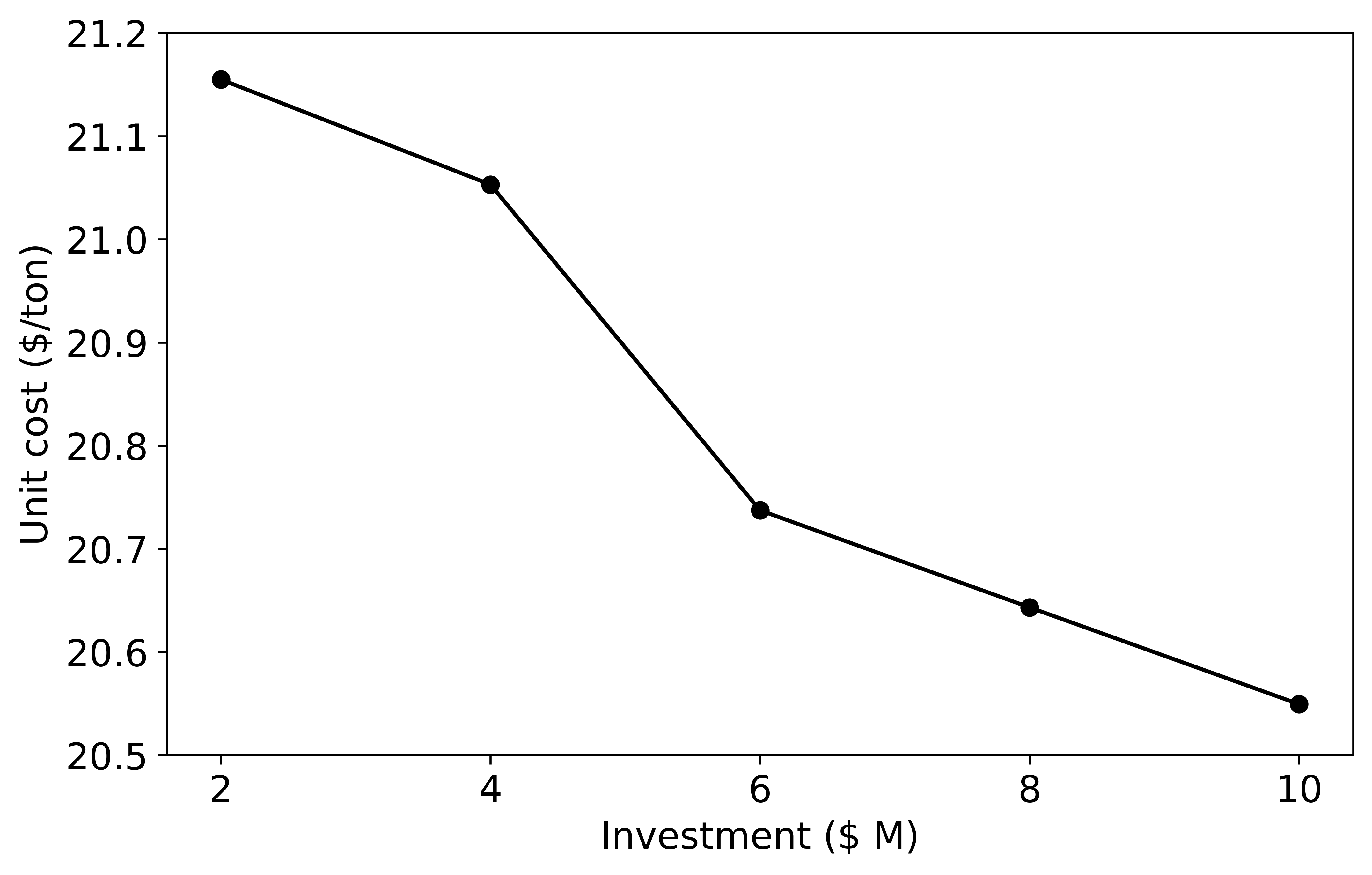}
      \caption{Unit Cost}
      \label{fig:UnitSupplyChainCost}
    \end{subfigure}
    \end{minipage}

    \caption{Investment vs. cost: (a) investment vs. total system cost, (b) investment vs. unit cost}
    \label{system_unit_cost}
\end{figure}

\begin{figure}[!htb]
    \centering

    \begin{minipage}[t]{.5\textwidth}
      \begin{subfigure}{\linewidth}
      \centering
      \includegraphics[width=1\linewidth]{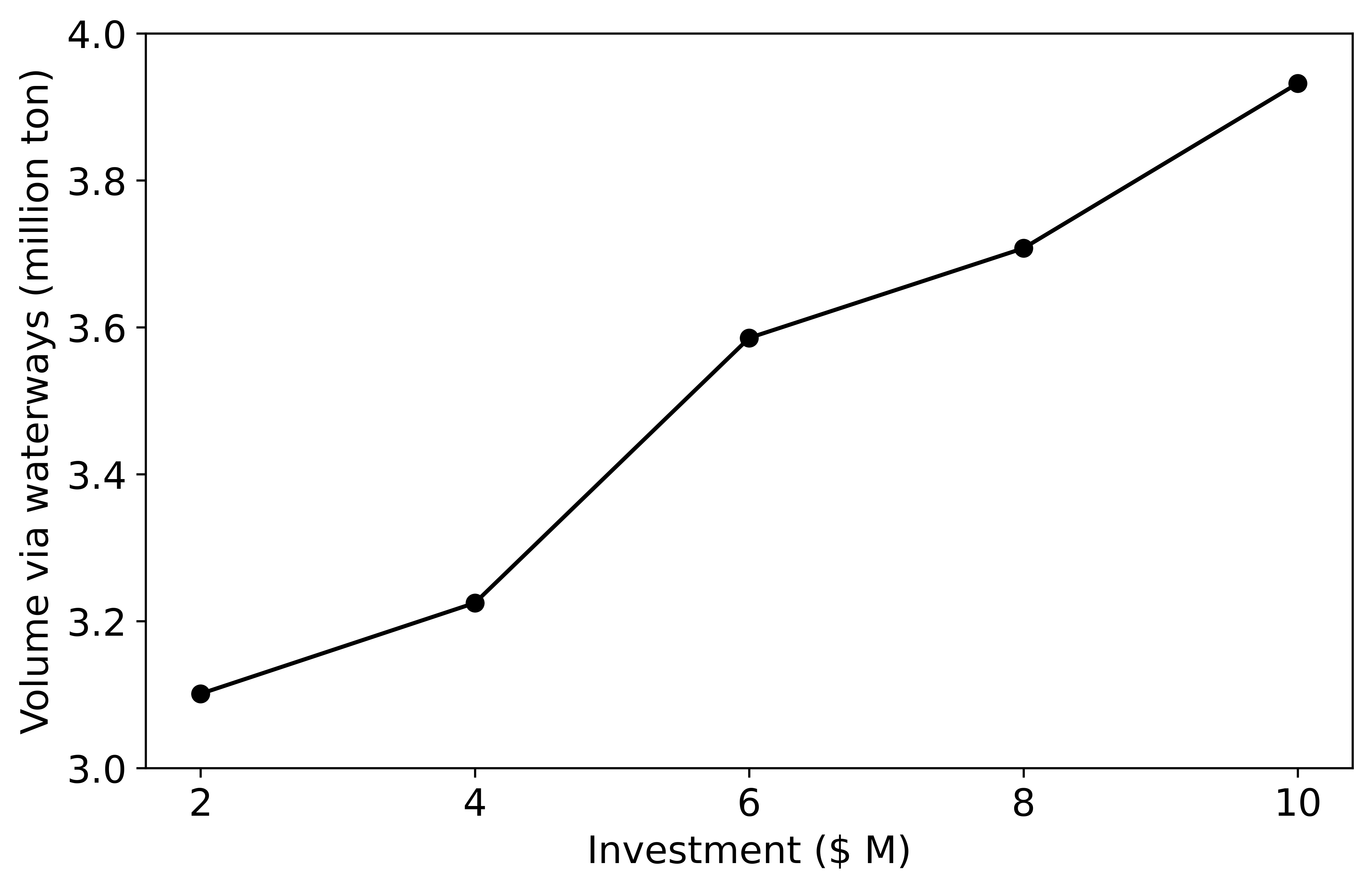}
      \caption{Volume}
      \label{fig:tonnage}
    \end{subfigure}
    \end{minipage}%
    \begin{minipage}[t]{.5\textwidth}
      \begin{subfigure}{\linewidth}
      \centering
      \includegraphics[width=1\linewidth]{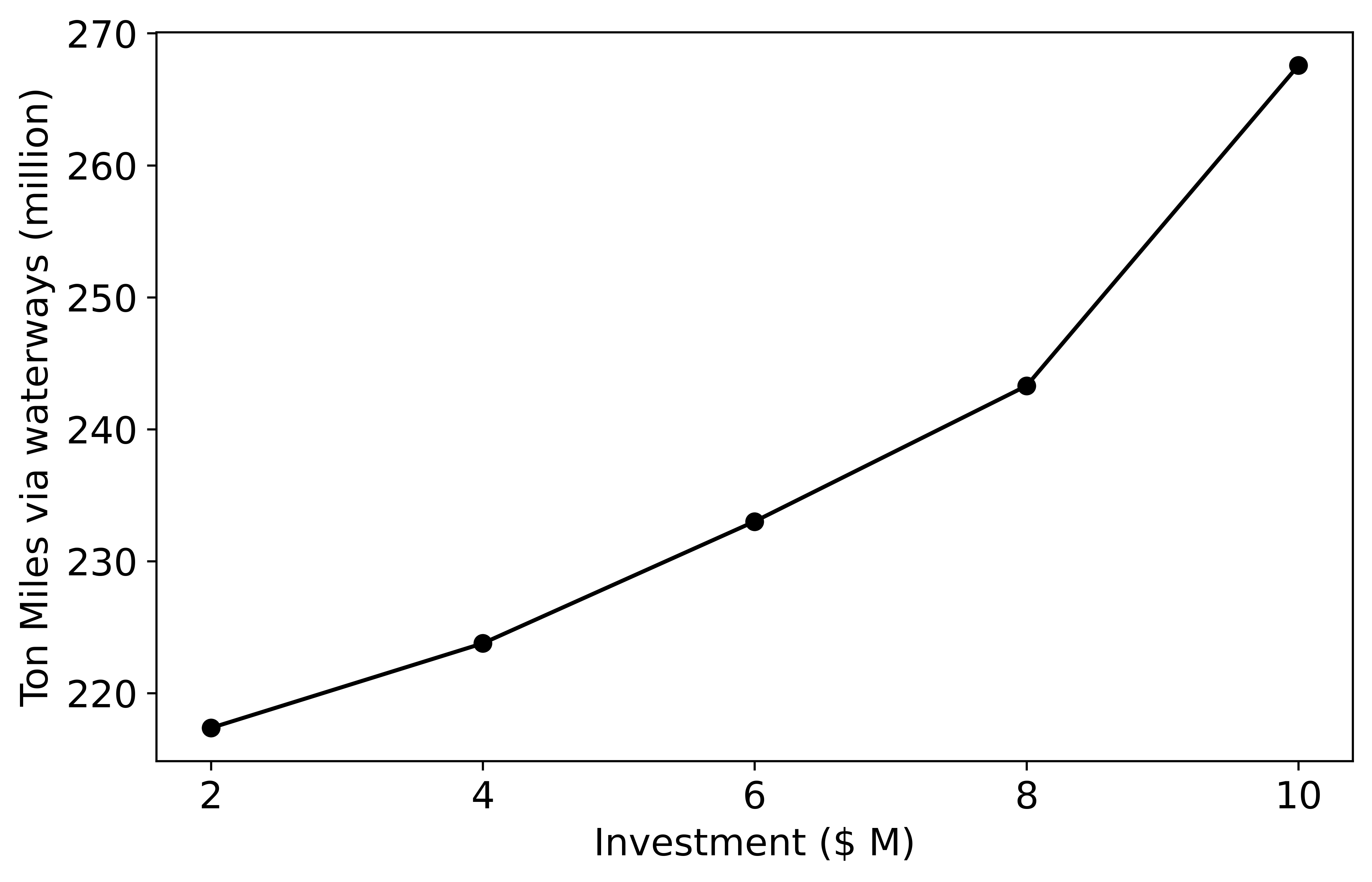}
      \caption{Ton-miles}
      \label{fig:tonmiles}
    \end{subfigure}
    \end{minipage}
        \begin{minipage}[t]{.5\textwidth}
      \begin{subfigure}{\linewidth}
      \centering
      \includegraphics[width=1\linewidth]{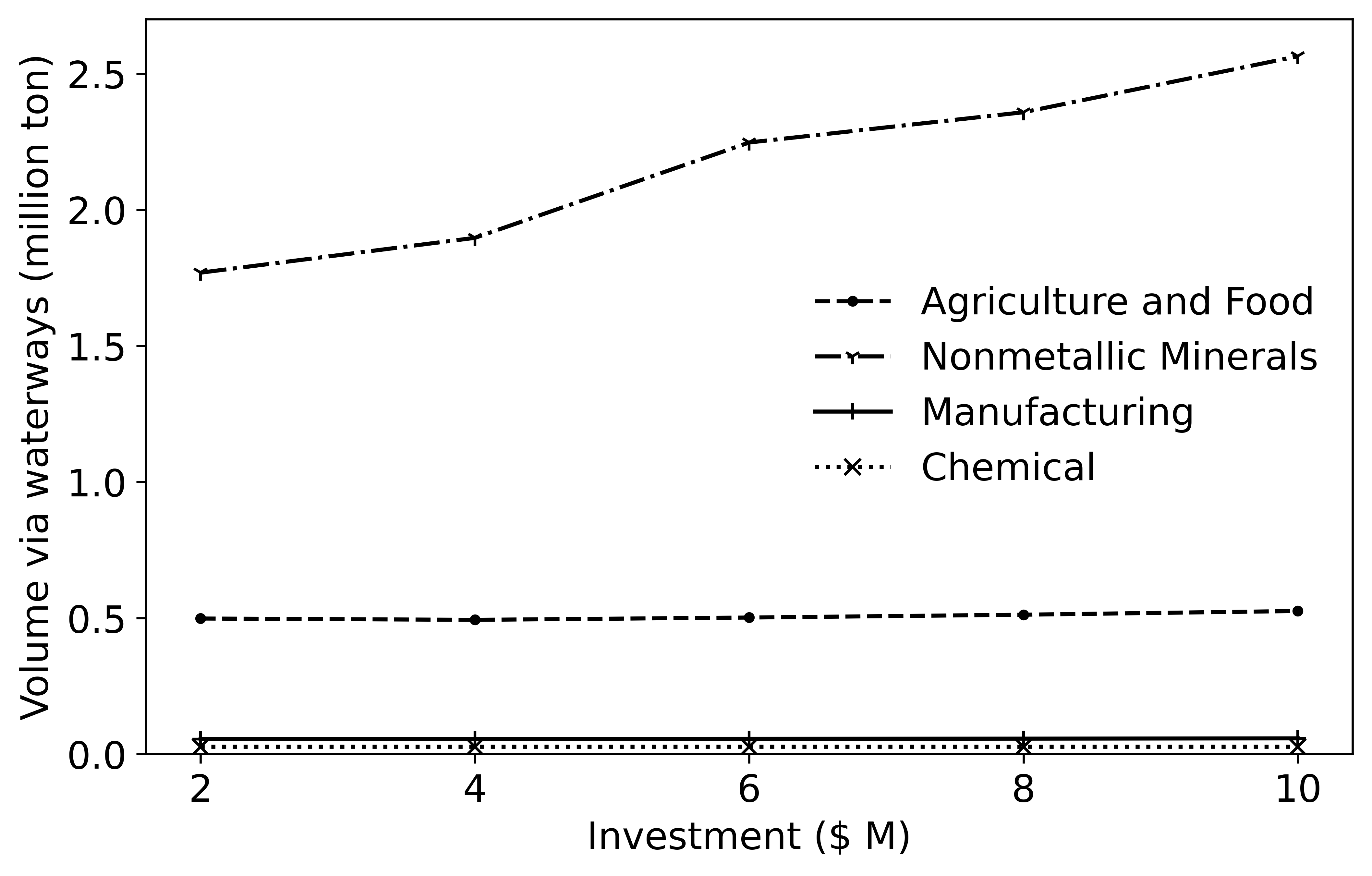}
      \caption{Commodity Volume}
      \label{fig:CommoditiesTonnage}
    \end{subfigure}
    \end{minipage}%
    \begin{minipage}[t]{.5\textwidth}
      \begin{subfigure}{\linewidth}
      \centering
      \includegraphics[width=1\linewidth]{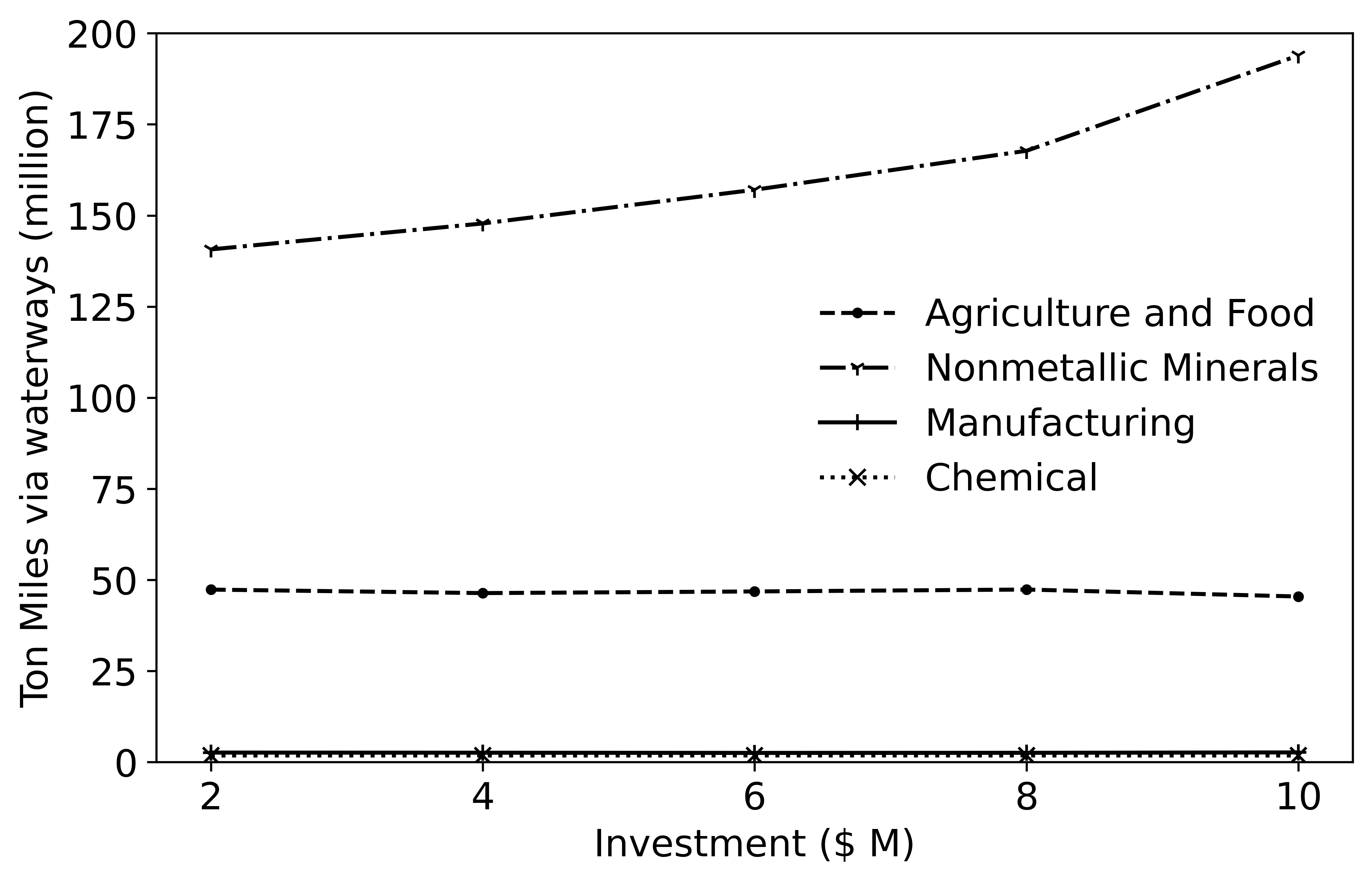}
      \caption{Commodity ton-miles}
      \label{fig:CommoditiesTonmiles}
    \end{subfigure}
    \end{minipage}
    \caption{Investment vs. volume and ton-miles via the waterways: (a) investment vs. total volume, (b) investment vs. total ton-miles, (c) investment vs. commodity volume, and (d) investment vs. commodity ton-miles}
    \label{system_supply_volume_tonmiles}
\end{figure}

We evaluate the impacts of varying infrastructure investments on system cost and volume of commodities moved via waterways. The total system cost, which includes transportation cost and investment cost, decreases from \$1.245 billion to \$1.225 billion (\$20 million decrease) as the investment increases from \$2 to \$6 million (\$4 million increase) (Figure \ref{system_unit_cost}). The total system cost drops to \$1.216 billion at an investment of \$10 million.  We see similar results for unit supply chain cost where the unit cost decreases from \$21.15/ton to \$20.74/ton (\$0.41/ton decrease) as the investment increases from \$2 million to \$6 million (\$4 million decrease) and eventually to \$20.55/ton for an investment of \$10 million  (Figure \ref{fig:UnitSupplyChainCost}).  This decrease in unit cost can be attributed to the increased volume of commodity shipped via waterways (Figure \ref{fig:tonnage}) which have the lowest transportation cost compared to truck and rail. This finding is consistent with previous research that reported economic benefits from expanding port infrastructure \cite{dekker2003economic}. From managerial perceptive, port authorities can use these findings as a part of benefit cost analysis to justify funding for inland waterway ports .\par

\begin{figure}[!htb]
    \centering
    \includegraphics[width = 1\linewidth]{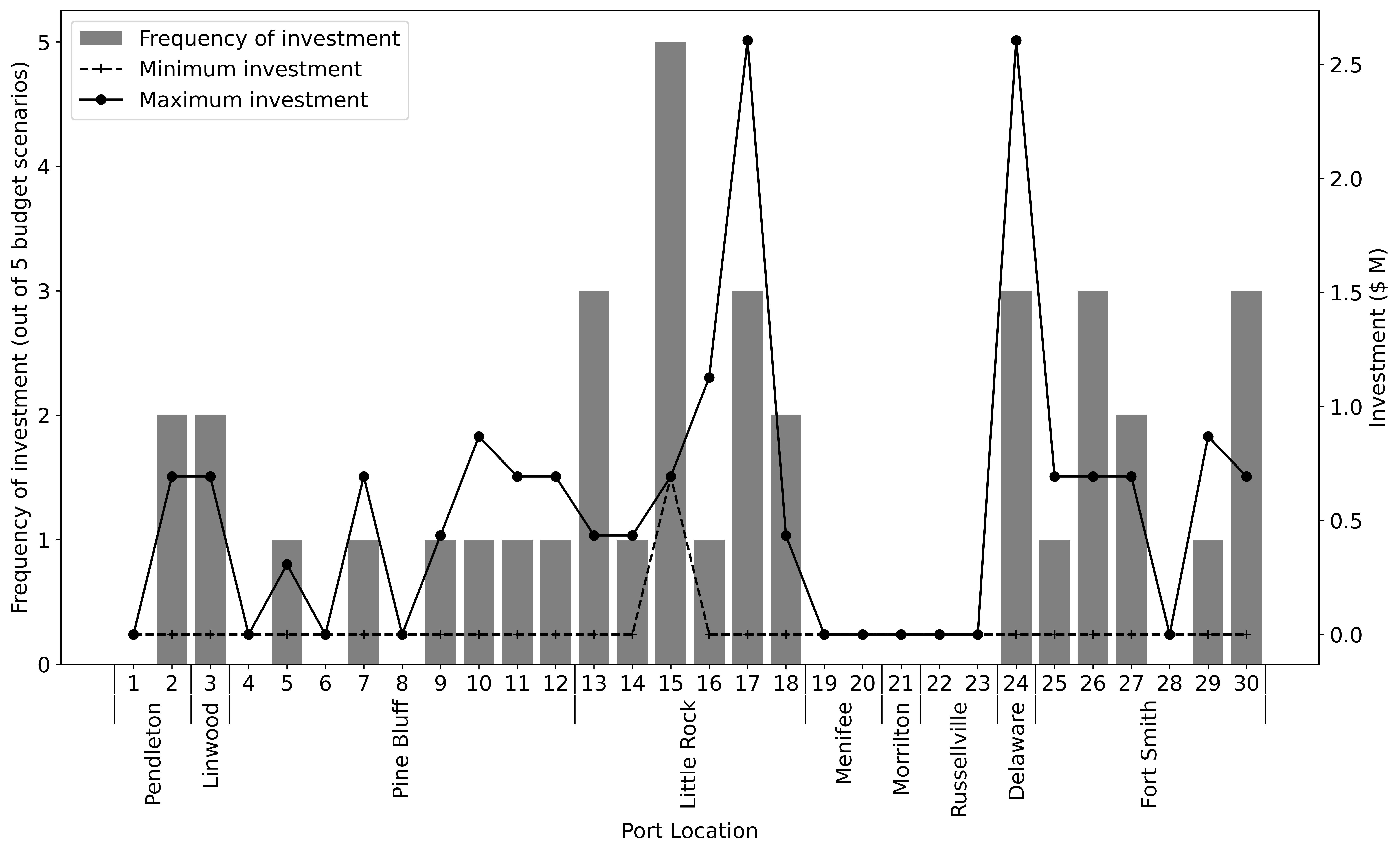}
    \caption{Port investment}
    \label{fig:portinvestment}
\end{figure}

We discuss the frequency and the range of investments at individual ports for  five different investment scenarios, \$2M, \$4M, \$6M, \$8M, and \$10M (Figure \ref{fig:portinvestment}). In most of these scenarios, the
model determines that investments should be made on ports located at Little Rock (central Arkansas), Delaware (midwestern section of the Arkansas river), and
Fort Smith (western most port on the Arkansas river). In every scenario considered, the model determines that investments in capacity expansion equipment should be made
at the port located at Little Rock (port 15 in Figure \ref{fig:portinvestment}). Across all scenarios, several ports repeatedly receive no investment e.g., Menifee, Morrilton, and Russellville. The findings suggest that the Little Rock port is an important hub of MKARNS, and this finding is corroborated by the prior investments of \$960,000 \$2 million  and \$3 million  received by the Little Rock Port Authority in 2012 \cite{us_economic_development_administration_us_2012}, in 2020 \cite{noauthor_dra_2020}, and in 2021 \cite{us_department_of_transportation_2021_2021}, respectively. From a managerial perspective, this result is crucial for port authorities to advocate for funds in the future. Similarly, port authorities can allocate those funds in the ports identified above for optimal movement of commodities through waterways in terms of unit supply chain cost as shown in Figure \ref{fig:UnitSupplyChainCost}.\par

The volume of commodities shipped via the waterways increases as the investment in port capacity increases (Figure (\ref{fig:tonnage})). A total of 2.74 million tons of commodities are shipped by waterways (6\% of the total shipments by the three modes) when the total investment is \$2 million. This volume increases to 3.47 million (8\% of the total shipments by all three modes) when the amount invested in infrastructure is \$10 million. We also report the change in ton-miles shipped by waterways for varying investment (Figure \ref{fig:tonmiles}). Ton-mile reflects both the volume (tons) and the distance (miles) shipped and is one of the most widely used measures of the physical volume of freight transportation services \cite{BTSTonMiles}. The waterway serves nearly 188 million ton-miles of freight (5\% of all three modes) at an investment of \$2 million and increases to 231 million ton-miles (6\% of all three modes) at an investment of \$10 million. The findings show that volume shipped via waterways and ton-miles increases with investment in port capacity expansions. Since barges have the least amount of carbon emission among other mode of transportation, from managerial perspective \cite{CarbonFootprintUSACE}, port stakeholders can use this finding to advocate for more funding for inland waterway port capacity expansion. \par

The 2-SOP model represents 11 commodity groups and models the shipment, supply, and demand for each commodity group uniquely while considering how port infrastructure can be shared among commodity groups, e.g., a warehouse can be used for manufacturing and primary metal commodities. To demonstrate this modeling contribution, we present the results for four commodity groups that dominate along the Arkansas section of the MKARNS \cite{asborno2020multicommodity, nachtmann2014economic}: nonmetallic minerals, agriculture and food, manufacturing, and chemicals. The volume of nonmetallic minerals shipped via waterways increases from 1.5 million (56\% of all commodities) to 2.2 million tons (64\% of all commodities)  (an increase of 0.8 million tons) as the investment increases from \$2 million to \$10 million (Figure \ref{fig:CommoditiesTonnage}). This represents a reduction of 1.4\% of non-metallic minerals shipped via truck and rail. Similarly, the nonmetallic mineral freight generated via waterways increases from 121 million ton-miles (64\% of all commodities) to 167 million ton-miles (72\% of all commodities) (an increase of 45 million ton-miles) when investments increase from \$2 million to \$10 million (Figure \ref{fig:CommoditiesTonmiles}). No significant impact of port infrastructure investment is seen for the other three commodity groups. This finding is crucial for the stakeholder, from a managerial perspective, to make decisions about their business approach. For example, since the findings show that the commodity transportation shifts away from trucks and trains and toward barges, companies related to nonmetallic minerals can use this information to make decisions on their business model regarding the number of trucks to acquire and the number of truck drivers to recruit.

\subsection{Evaluation of stochastic solutions}
To demonstrate the benefit of developing the model with stochastic elements, rather than a deterministic model, we calculate the Value of Stochastic Solution (VSS) as the difference between the objective function value of the stochastic solution and the expected value solution. In addition, we calculate the Expected Value of Perfect Information to quantify the value of perfect information required to predict future with certainty. EVPI is the difference between the objective function value of the stochastic solution and the wait-and-see solutions (WSS). \par
\begin{table}[]
\centering
\caption{Comparison of stochastic and deterministic solutions}
\label{tab:performanceevaluation}
\begin{tabular}{ll}
\hline
Strategies                  & Value (\$M)  \\ \hline
Stochastic   programming    & 1,225        \\
\multicolumn{2}{l}{Wait and see solutions} \\
Scenario-1                  & 1,433        \\
Scenario-2                  & 626          \\
Scenario-3                  & 680          \\
Scenario-4                  & 538          \\
Scenario-5                  & 1,367        \\
Scenario-6                  & 2,336        \\
Scenario-7                  & 976          \\
Scenario-8                  & 286          \\
Scenario-9                  & 1,455        \\
Scenario-10                 & 1,798        \\
Expected value   solution   & 1,246        \\ \hline
\end{tabular}
\end{table}

The expected value solution is calculated in two steps. First, we solve (WSN) assuming a
single scenario, represented by the expected values of commodity demand. Next, we fix the
values of the first stage decisions using this solution, and resolve (WSN) to obtain the expected value solution of \$1,246M. Hence, VSS = \$1,246M - \$1,225M = \$21M. Therefore, \$21M per year can be saved by solving the stochastic model compared to the corresponding deterministic model.

The following is the approach we use to calculate EVPI. Let us assume that we know
exactly what scenario is realized in the future. Then, we can solve (WSN) for this particular scenario. This is the “wait-and-see solution" (WSS). Table \ref{tab:performanceevaluation} summarizes the WSS of each scenario. We calculate the expected value of WSS to be \$1,221M (\$1,433M*0.065 + \$626M*0.07 + \$680M*0.075 + \$538M*0.08 + \$1,367M*0.085 + \$2,336M*0.09 + \$976M*0.095 + \$286M*0.1 + \$1,455M*0.20 + \$1,798M*0.14). Recall the multipliers in the above calculation correspond to the likelihood of each scenario occurring. Therefore, EVPI is \$4M (\$1,225M - \$1,221M ). Therefore, the price we should be willing to pay to correctly predict future commodity demand realizations should not exceed \$4M.  \par 


\section{Conclusions}
In this paper, we develop a model to guide strategic investments in inland waterway port infrastructure investments given uncertainty in commodity demand. This study proposes a 2-SOP model that seeks to minimize the total of port infrastructure investment costs and the expected transportation costs. We implement a Benders decomposition algorithm to solve the model, and accelerate this algorithm using Knapsack inequalities and Pareto-optimal cuts. The computational analysis reveals that Benders with Knapsack inequalities and Pareto-optimal cuts outperforms Gurobi and the traditional Benders algorithm for large-sized problems.
 \par
We apply the two-stage stochastic optimization model to the Arkansas section of the MKARNS. The model results show that while the total system cost (transportation plus investment costs) decreases with increasing investment, the rate of decrease in system cost is convex in nature, i.e., the rate of change decreases with the dollar amount invested in port capacity expansion. Our model shows that commodity
volume and, as expected, the percent of that volume that moves via waterways (in ton-miles) increases with increasing investment in port infrastructure. The model captures individual commodity movements in terms of tons and ton-miles shipped by transportation mode.  Results show that among all commodities, nonmetallic minerals experience the largest fluctuation in the tonnage and ton-miles shipped as a consequence of changing investment amounts. Furthermore, since the model estimates investments and commodity throughput at individual ports, we are able to identify a cluster of ports (Little Rock, Fort Smith) that should receive investment in port capacity under any investment scenario. \par
Finally, to demonstrate the value of a stochastic formulation over a deterministic approach, we calculate the value of the stochastic solution, VSS. VSS shows that a failure to use stochastic model to capture variations in commodity demand, could cost up to \$21 M per year. We also calculate the expected value of perfect information (EVPI). EVPI indicates the price we should be willing to pay for correctly predicting future realizations of commodity demand, should be no more than \$4M. \par
This research opens several avenues for research to explore in the future. The current model allows for analysis of the supply chain within a state, and hence in the future, the model will be expanded to incorporate multiple states. Furthermore, in our current study, there are no port disruption scenarios considered, and since closure of port due to human-induced or / and natural causes can impact port operation, the uncertainty of port disruption will be incorporated in a future study. While our current work used knapsack inequalities and Pareto-optimal cuts to accelerate the Benders decomposition algorithm, our future work will explore additional enhancement techniques such as maximum density cut generation and Benders-type heuristics as mentioned by \citep{rahmaniani2017benders}. 

\bibliographystyle{chicago}
\bibliography{trb_template}
 
\appendix
\section{Appendix A}
\label{appendix:dual}
\begin{subequations}
\small
\begin{align}
\begin{split}
\label{equ:dualsecondstage}(DS-WSN(s)): & \max \sum_{j \in J}\sum_{c\in C}\sum_{p \in P}\sum_{s \in S} (\nu_{jcps}q_{jcps} + \xi_{jcps}d_{jcps}) + 
\sum_{i \in I}\sum_{e\in E}\sum_{p \in P}\sum_{s \in S}\pi_{ipes}m_{e}(n_{ie}+Z_{ie})+\\
&\hspace{0.18 in}\sum_{i \in I}\sum_{f\in F}\sum_{p \in P}\sum_{s \in S}\sigma_{ipfs}l_f(k_{if}+Y_{if})
\end{split}
\end{align}
Subject to:
\small
\begin{align}
    \label{cons:DS_xt}&\nu_{jcps}+\upsilon_{icps}+\Lambda_{ec}\pi_{ipes}\leq \alpha ^t_{ij},\hspace{0.1in} \forall i \in I^{\prime}, j \in J^{\prime}, c \in C, e \in E, p \in P, s \in S,\beta_{ec} = 1 \\
    \label{cons:DS_xr}&\nu_{jcps}+\upsilon_{icps}+\Lambda_{ec}\pi_{ipes}\leq \alpha ^r_{ij},\hspace{0.1in}\forall i \in I^{\prime}, j \in J^{\prime}, c \in C, e \in E, p \in P, s \in S,\beta_{ec} = 1, \gamma_j = 1, \delta_i = 1\\
    \label{cons:DS_u1}&\upsilon_{icp+1s}-\upsilon_{icps}+\zeta_{fc}\sigma_{ipfs}\leq h_{c},\hspace{0.1in}\forall i \in I^{\prime}, c \in C,f \in F, p \in 1..P-1, s \in S, \Gamma_{fc}=1\\
    \label{cons:DS_u2}&-\upsilon_{icps}+\zeta_{fc}\sigma_{ipfs}\leq h_{c},\hspace{0.1in}\forall i \in I^{\prime}, c \in C,f \in F, p \in |P|, s \in S, \Gamma_{fc}=1\\
    \label{cons:DS_w}&\chi_{kcps}-\upsilon_{icps}+\Lambda_{ec}(\pi_{ipes}+\pi_{kpes})\leq a_{ik},\hspace{0.1in}\forall i \in I^{\prime},k \in I^{{\prime}{\prime}}, c \in C,e \in E, p \in P, s \in S,\beta_{ec} = 1\\
    \label{cons:DS_v1}&\chi_{icp+1s}-\chi_{icps}+\zeta_{fc}\sigma_{ipfs}\leq h_{c},\hspace{0.1in}\forall i \in I^{{\prime}{\prime}}, c \in C,f \in F, p \in 1..P-1, s \in S, \Gamma_{fc}=1\\
    \label{cons:DS_v2}&-\chi_{icps}+\zeta_{fc}\sigma_{ipfs}\leq h_{c},\hspace{0.1in}\forall i \in I^{{\prime}{\prime}}, c \in C,f \in F, p \in |P|, s \in S, \Gamma_{fc}=1\\
    \label{cons:DS_rt}&\xi_{jcps}-\chi_{icps}+\Lambda_{ec}\pi_{ipes}\leq \alpha ^t_{ij},\hspace{0.1in}\forall i \in I^{{\prime}{\prime}},j \in J^{{\prime}{\prime}}, c \in C,e \in E, p \in P, s \in S,\beta_{ec} = 1\\
    \label{cons:DS_rr}&\xi_{jcps}-\chi_{icps}+\Lambda_{ec}\pi_{ipes}\leq \alpha ^r_{ij},\hspace{0.1in}\forall i \in I^{{\prime}{\prime}},j \in J^{{\prime}{\prime}}, c \in C,e \in E, p \in P, s \in S,\beta_{ec} = 1, \gamma_j = 1, \delta_i = 1\\
    \label{cons:DS_ot}&\nu_{jcps}+\xi_{jcps}\leq l^t_{jm},\hspace{0.1in}\forall j \in J^{\prime},m \in J^{{\prime}{\prime}},  c \in C, p \in P, s \in S\\
    \label{cons:DS_or}&\nu_{jcps}+\xi_{jcps}\leq l^r_{jm},\hspace{0.1in} \forall j \in J^{\prime},m \in J^{{\prime}{\prime}}, c \in C, p \in P, s \in S, \gamma_j = 1, \gamma_{j^{\prime}} = 1\\
    \label{cons:DS_e}&\xi_{jcps} \leq \mu,\hspace{0.1in} \forall j \in J, c \in C, p \in P, s \in S\\
    \label{cons:DS_bounds}&\nu_{jcps}, \pi_{ipes}, \sigma_{ipfs} \in R_{-}\\
    \label{cons:DS_free}&\xi_{jcps},\upsilon_{icps}, \chi_{icps} \in free
 \end{align}
 \end{subequations}
 \clearpage
\section{Appendix B}
\label{appendix:a}
\begin{table}[!ht]
    \centering
    \caption{Transportation cost for truck, rail and barge}
    \begin{tabular}{ccc}
    \hline
        Transportation Mode & Parameters & Value  \\ \hline
        Truck & $\alpha^{t}_{ij}$, $l^{t}_{jm}$ & \$0.185/ton-mile  \\ 
        Rail & $\alpha^{r}_{ij}$, $l^{r}_{jm}$ & \$22.65/ton+\$0.033/ton-mile  \\ 
        Barge & $a^{r}_{ik}$ & \$0.0089/ton-mile  \\ \hline
    \end{tabular}
    \label{tab:transportationcost}
\end{table}

\begin{table}[!ht]
    \centering
    \caption{Equipment cost}
    \begin{tabular}{ccc}
    \hline
        Equipment & Cost (\$) & Specification  \\ \hline
        Conveyor &  18,723  & 36 inch  \\ 
        Crane &  300,000  & 65 ton  \\ 
        Hopper &  18,723  & 25 ton  \\ 
        Forklift &  96,738  &   \\ \hline
    \end{tabular}
    \label{tab:equipmentcost}
\end{table}

\begin{table}[!ht]
    \centering
     
    \caption{Storage facility cost}
    \begin{tabular}{ccc}
    \hline
        Storage facility & Cost (\$) & Specification  \\ \hline
        Grain elevator &  227,866  & 650,000 bushels  \\ 
        Unpaved storage &  692,769  & 547,000 $ft^2$  \\
        Paved storage &  307,065  & 144,000 $ft^2$  \\ 
        Warehouse &  5,663,854  & 77,000 $ft^2$  \\ 
        Chemical/Petroleum storage tank &  1,109,090  & 125,000 barrels  \\ \hline
    \end{tabular}
    \label{tab:storagefacilitycost}
\end{table}

\begin{table}[!t]
    \centering
    \caption{Equipment processing capacity (ton/month)}
    \begin{tabular}{p{1cm}*{5}{p{2cm}}}
    \hline
        \multirow{4}{*}{Port} &  Crane, Conveyor, Hopper, Forklift  & Crane, Forklift & Petroleum tank & Chemical tank  \\ \hline
        1 & 32,400 & 0 & 0 & 7,624  \\ 
        2 & 34,425 & 0 & 0 & 0  \\ 
        3 & 30,375 & 0 & 0 & 0  \\ 
        4 & 30,000 & 0 & 0 & 0  \\ 
        5 & 0 & 150,000 & 0 & 7,624  \\ 
        6 & 0 & 0 & 0 & 0  \\ 
        7 & 30,000 & 0 & 0 & 0  \\ 
        8 & 50,250 & 9,300 & 0 & 0  \\ 
        9 & 30,000 & 0 & 0 & 0  \\ 
        10 & 0 & 210,000 & 0 & 0  \\ 
        11 & 30,000 & 0 & 0 & 0  \\ 
        12 & 30,375 & 0 & 0 & 0  \\ 
        13 & 38,700 & 200,700 & 0 & 0  \\ 
        14 & 0 & 0 & 38,120 & 0  \\ 
        15 & 129,300 & 21,600 & 7,624 & 0  \\ 
        16 & 0 & 0 & 0 & 22,872  \\ 
        17 & 105,000 & 0 & 0 & 0  \\ 
        18 & 0 & 58,500 & 0 & 0  \\ 
        19 & 26,250 & 0 & 0 & 0  \\ 
        20 & 52,500 & 0 & 0 & 0  \\ 
        21 & 30,000 & 0 & 0 & 0  \\ 
        22 & 26,250 & 0 & 0 & 0  \\ 
        23 & 4,050 & 0 & 0 & 0  \\ 
        24 & 52,500 & 0 & 0 & 0  \\ 
        25 & 15,000 & 0 & 0 & 0  \\ 
        26 & 0 & 90,000 & 0 & 0  \\ 
        27 & 76,650 & 30,000 & 0 & 0  \\ 
        28 & 20,250 & 0 & 0 & 0  \\ 
        29 & 34,425 & 0 & 0 & 0  \\ 
        30 & 105,000 & 0 & 0 & 0  \\ \hline
    \end{tabular}
    \label{tab:processingcapacity}
\end{table}

\begin{table}[!t]
    \centering
    \caption{Storage facility capacity (ton)}
    \begin{tabular}{p{1cm} *{7}{p{1.8cm}}}
    \hline
        \multirow{2}{*}{Port} & Grain Elevator & Unpaved Storage & Paved Storage & Warehouse & Chemical storage tank & Petroleum storage tank  \\ \hline
        1 & 118,800 & 18,687 & 0 & 4,182 & 0 & 3,600  \\ 
        2 & 15,984 & 0 & 0 & 0 & 0 & 0  \\ 
        3 & 61,992 & 0 & 0 & 0 & 0 & 0  \\ 
        4 & 11,556 & 0 & 0 & 0 & 0 & 0  \\ 
        5 & 0 & 0 & 0 & 15,410 & 0 & 0  \\ 
        6 & 0 & 0 & 0 & 0 & 0 & 26,250  \\ 
        7 & 11,214 & 0 & 0 & 0 & 0 & 0  \\ 
        8 & 324 & 0 & 0 & 48,956 & 0 & 0  \\ 
        9 & 0 & 176,380 & 0 & 0 & 0 & 0  \\ 
        10 & 0 & 115,352 & 0 & 3,679 & 0 & 0  \\ 
        11 & 0 & 0 & 0 & 4,594 & 0 & 0  \\ 
        12 & 113,400 & 0 & 0 & 0 & 0 & 0  \\ 
        13 & 0 & 143,749 & 191,602 & 5,906 & 0 & 0  \\ 
        14 & 0 & 0 & 0 & 0 & 29,700 & 0  \\ 
        15 & 56,700 & 0 & 5,748,048 & 10,731 & 7,950 & 0  \\ 
        16 & 0 & 0 & 0 & 0 & 0 & 27,300  \\ 
        17 & 0 & 261,766 & 0 & 0 & 0 & 0  \\ 
        18 & 0 & 9,793 & 45,646 & 0 & 0 & 0  \\ 
        19 & 0 & 48,243 & 0 & 0 & 0 & 0  \\ 
        20 & 0 & 50,614 & 188,185 & 0 & 0 & 0  \\ 
        21 & 13,500 & 0 & 0 & 1,254 & 0 & 0  \\ 
        22 & 0 & 0 & 316,559 & 0 & 0 & 0  \\ 
        23 & 17,550 & 0 & 0 & 10,073 & 0 & 0  \\ 
        24 & 0 & 1,069,815 & 0 & 0 & 0 & 0  \\ 
        25 & 0 & 0 & 1,322,985 & 4,534 & 0 & 0  \\ 
        26 & 0 & 0 & 0 & 10,047 & 0 & 0  \\ 
        27 & 0 & 0 & 492,369 & 13,922 & 0 & 0  \\ 
        28 & 22,950 & 0 & 0 & 0 & 0 & 0  \\ 
        29 & 0 & 125,172 & 0 & 0 & 0 & 0  \\ 
        30 & 0 & 0 & 0 & 25,988 & 0 & 0  \\ \hline
    \end{tabular}
    \label{tab:storagecapacity}
\end{table}
\end{document}